\documentclass[11pt]{article}
\usepackage{amsfonts}
\usepackage{amscd}
\usepackage{amsmath}
\usepackage{epsfig}

\newlength{\dinwidth}
\newlength{\dinmargin}
\newtheorem{theorem}{Theorem}
\newtheorem{proposition}{Proposition}
\newtheorem{corollary}{Corollary}
\newtheorem{remark}{Remark}
\newtheorem{lemma}{Lemma}

\def\be{\begin{equation}}
\def\ee{\end{equation}}
\def\ben{\begin{displaymath}}
\def\een{\end{displaymath}}
\def\baa{\begin{eqnarray}}
\def\eaa{\end{eqnarray}}

\def\ba{\begin{array}}
\def\ea{\end{array}}
\makeatletter
\@addtoreset{equation}{section}
\makeatother

\def\cdiff{{\cal C}}
\def\Wcal{{\cal W}}
\def\CP1{{\bf CP}^1}

\def\C{{\bf C}}

\def\X{{\cal X}}
\def\a{\alpha}

\def\b{\beta}

\def\l{z}

\def\P{\wp}

\def\la{\label}

\def\L{{\cal X}}

\def\p{\partial}
\def\B{{\bf B}}

\def\P{P}

\begin{document}
\title{On the asymptotics of determinant of Laplacian at the principal boundary of the principal stratum of the moduli space of Abelian differentials}

\author{A. Kokotov\footnote{e-mail: alexey@mathstat.concordia.ca} }

\maketitle

\begin{center}
Department of Mathematics and Statistics, Concordia University\\
1455 de Maisonneuve Blvd. West \\
Montreal, Quebec H3G 1M8 Canada
\end{center}


Let $\X$ be a translation surface of genus $g>1$ with $2g-2$ conical points of angle $4\pi$ and let $\gamma$, $\gamma'$ be two  homologous saddle connections of length $s$ joining two conical points of $\X$ and bounding two surfaces $S^+$ and $S^-$ with boundaries $\partial S^+=\gamma-\gamma'$ and $\partial S^-=\gamma'-\gamma$.
Gluing the opposite sides of the boundary of each surface $S^+$, $S^-$ one gets two (closed) translation surfaces  $\X^+$, $\X^-$ of genera $g^+$, $g^-$; $g^++g^-=g$. Let $\Delta$,  $\Delta^+$ and $\Delta^-$ be the Friedrichs extensions of the Laplacians corresponding to the (flat conical) metrics on $\X$, $\X^+$ and $\X^-$ respectively. We study the asymptotical behavior of the (modified, i. e. with zero modes excluded) zeta-regularized determinant  ${\rm det}^*\, \Delta$ as  $\gamma$ and $\gamma'$ shrink. We find the asymptotics $${\rm det}^*\,\Delta\sim \kappa s^{1/2}\frac{{\rm Area}\,(\X)}{{\rm Area}\,(\X^+){\rm Area}\,(\X^-)}\,{\rm det}^*\,\Delta^+{\rm det}^*\,\Delta^-$$
as $s\to 0$; here $\kappa$ is a certain absolute constant admitting an explicit expression through spectral characteristics of some model operators.
We use the obtained result to fix an undetermined constant in the explicit formula for ${\rm det}^*\, \Delta$ found in \cite{DG}.

\section{Introduction}

Let $H_g(1, \dots, 1)$ (one has here $2g-2$ units) be the principal stratum of the moduli space of Abelian differentials over compact Riemann surfaces of genus g. One defines $H_g(1, \dots, 1)$ as the moduli space of pairs $(\X, \omega)$, where $\X$ is a compact Riemann surface of genus $g$ and $\omega$ is a holomorphic one-form (an Abelian differential) on $\X$  with $2g-2$ zeros of multiplicity one. It is known (\cite{KZ}) that $H_g(1, \dots, 1)$ is a connected complex orbifold of (complex) dimension $4g-3$.

  Let a pair $(\X, \omega)$ belong to $H_g(1, \dots, 1)$. The holomorphic differential $\omega$ defines the conformal flat conical metric $|\omega|^2$ on  $\X$, this metric has conical points of angle $4\pi$ at the zeros of $\omega$ and trivial monodromy along any closed loop in $\X\setminus \{{\rm conical \ points}\}$. Thus, the 2-d manifold $\X$ equipped with metric $|\omega|^2$ becomes a so-called {\it translation surface}. It should be noted that any translation surface (a compact 2-manifold with flat conical metric having trivial holonomy) can be obtained as a pair $(\X, |\phi|^2)$, where $\X$ is a compact Riemann surface and $\phi$ is an Abelian differential on $\X$ (in general, with zeros of arbitrary multiplicities).

To the metric $|\omega|^2$ one can associate the Laplace operator $\Delta^{|\omega|^2}$ (often denoted below simply by $\Delta$; we assume $\Delta$ be a nonnegative operator, i. e. one attaches minus to the usual definition of the Laplace operator) with domain $C^\infty_c(\X\setminus \{{\rm conical \ points}\})$. The Friedrichs extension of $\Delta$ (from now on the notation $\Delta$ refers  only to this self-adjoint operator in the Hilbert space $L_2(\X, |\omega|^2))$ is known to have discrete spectrum $0=\lambda_0<\lambda_1\leq \lambda_2\leq \dots$ of finite multiplicity. The operator zeta-function, defined for $\Re t>1$ as
$$\zeta_\Delta(t)=\sum_{\j=1}^\infty\lambda_j^{-t}$$
admits analytic continuation to ${\mathbb C}$ as a meromorphic function with the only pole $t=1$. The (modified) zeta-regularized determinant of the operator $\Delta$ is defined via the relation $\log {\rm det}^*\,\Delta=-\zeta'_\Delta(0)$.

If $\X$ is an elliptic curve the Abelian differentials on $X$ have no zeros; the moduli space of Abelian differentials on Riemann surfaces of genus one is denoted by $H_1(\emptyset)$. Introduce the real-valued function ${\cal F}_1$ on $H_1(\emptyset)$ via
$$H_1(\emptyset)\ni (\X, \omega)\mapsto {\cal F}_1(\X, \omega)={\rm det}^*\,\Delta^{|\omega|^2}\,.$$
In genus one the spectrum of the operator $\Delta$ is known explicitly and the direct calculation of the value $\zeta'_\Delta(0)$  (which essentially reduces to making use of the first Kronecker limit formula) leads to the following expression (found in \cite{RS}; see also \cite{Polch})
\begin{equation}\label{Ray}{\cal F}_1(\X, \omega)=
4{\Im(B/A){\rm Area}(\X, |\omega|^2)}|\eta(B/A)|^4\,,
\end{equation}
where $A=\oint_a\omega$, $B=\oint_b\omega$ with $\{a, b\}$ being a canonical basis of cycles on $\X$; ${\rm Area}(\X, |w|^2)=\Im(A\bar{B})$,  and $\eta$ is the Dedekind eta-function
$$\eta(\sigma)=\exp\left(\frac{\pi i \sigma}{12}\right)\prod_{n\in {\mathbb N}}\big(1-\exp(2\pi i n\sigma)\big)\,.$$
In \cite{DG} this classical result was generalized to the case of an arbitrary genus and an explicit expression for the function
$$H_g(1, \dots, 1)\ni (\X, \omega)\mapsto {\cal F}_g(\X, \omega)={\rm det}^*\,\Delta^{|\omega|^2}$$
was found. To formulate this result we need to introduce some auxiliary objects.
Let $\{a_\a, b_\a\}_{\a=1, \dots, g}$ be a canonical basis of cycles on $\X$. Denote by $\widehat \X$ a fundamental polygon obtained
via cutting the surface $\X$ along a system of $2g$ loops starting at some chosen point of $\X$ and homologous to the basic cycles.

 Introduce the basis of normalized Abelian differentials $\{v_\alpha\}$ on $\X$, the matrix of corresponding $b$-periods $\B=||\oint_{\b_\alpha} v_\beta||$ and the vector of Riemann constants:
\be
K^P_\alpha=\frac{1}{2}+\frac{1}{2}\B_{\a\a}-\sum_{\b=1,
\b\neq\a}^g\oint_{a_\b}\left(v_\beta\int_P^xv_\a\right)\;,
\la{rimco}\ee
where the interior integral is taken along a path which does not
intersect $\partial\widehat \X$. Let $E(P, Q)$ be the Schottky-Klein prime form (see \cite{Fay}).

As in \cite{Fay92} introduce
\begin{itemize}
\item

 the following holomorphic multi-valued $(g/2, -g/2)$-differential  $\sigma(P,Q)$:
\begin{equation}\label{sigma}
\sigma(P, Q)=\exp\left\{
-\sum_{\a=1}^g\oint_{a_\a}v_\a(R)\log\frac{E(R, P)}{E(R, Q)}
\right\}\,;
\end{equation}
 the right-hand side of
(\ref{sigma}) is a non-vanishing holomorphic $g/2$-differential on
$\widehat \X$ with respect to $P$ and a non-vanishing holomorphic
$(-g/2)$-differential with respect to $Q$;
\item
 the following holomorphic multivalued
$g(1-g)/2$-differential on $\X$:
\begin{equation}\label{c}
\cdiff(P)=\frac{1}{\Wcal[v_1, \dots, v_g](P)}\sum_{\a_1, \dots,
\a_g=1}^g
\frac{\partial^g\Theta(K^P)}{\p z_{\a_1}\dots \p z_{\a_g}}
v_{\a_1}\dots v_{\a_g}(P)\;,
\end{equation}
where
\be
\Wcal(P):= {\rm \det}_{1\leq \a, \b\leq g}||v_\b^{(\a-1)}(P)||
\la{Wronks}
\ee
is the Wronskian determinant of holomorphic differentials
 at the point $P$.

\end{itemize}

Let $(\omega)=\sum_{k=1}^{2g-2}P_k$ be the divisor of the holomorphic differential $\omega$,
  denote by ${\cal
A}_P(\cdot)$ the Abel map with the base point $P$. Then one has the relation
\begin{equation}\label{rdef}{\cal A}((\omega))+2K^P+{\mathbb B}{\bf r}+{\bf
q}=0\end{equation}
 with some integer vectors ${\bf r}$ and ${\bf q}$.
 Let us emphasize that vectors ${\bf r}$, ${\bf q}$ as well as the prime form and the differentials $C$ and $\sigma$ depend on the choice of the
 fundamental polygon $\widehat \X$.

Now we are able to formulate the result from \cite{DG}. One has the relation
\begin{equation}\label{main}{\cal F}_g(\X, \omega)=\delta_g{\rm det}\,\Im \B\,{\rm Area}(\X, |\omega|^2)|\tau_g(\X, \omega, \{a_\a, b_\a\})|^2\,\end{equation}
where $\delta_g$ is a constant depending only on genus $g$ and $\tau_g(\X, \omega, \{a_\a, b_\a\})$ is defined up to a unitary multiplicative factor (and not a choice of the fundamental polygon!) by the formula
\begin{equation}\label{deftau}
\tau^{-6}_g(\X, \omega, \{a_\a, b_\a\})=e^{2\pi i<{\bf r},
K^P>}C^{-4}(P)\prod_{k=1}^{2g-2}\sigma(P_k, P)\left\{E(P_k, P)
\right\}^{(g-1)}\,.
\end{equation}
Here $P$ is an arbitrary point of $\X$ and the integer vector ${\bf r}$ is
defined by (\ref{rdef}),  the values of the prime form and $\sigma$ at the zeros $P_k$ of the differential $\omega$ are calculated in the local parameter $x_k(Q)=\sqrt{\int_{P_k}^Q\omega}\,$, the values of the prime form and $\sigma$ at the point $P$ are taken in the local parameter $z(Q)=\int^Q\omega$; the expression (\ref{deftau}) is independent of the choice of $P$.

\begin{remark}{\rm In case $g=1$, using (\ref{Ray}), the formula
 $$C(P)=2\pi i\eta^3(B/A)e^{-\pi i\frac{B}{4A}}$$
from (\cite{Fay92}, p. 21) and the relation $K^P=\frac{1}{2}+\frac{B}{2A}$
(implying ${\bf r}=-1$ in (\ref{rdef})) together with (\ref{main}) and (\ref{deftau}), one gets the relation
\begin{equation}\label{delta1}
\delta_1=\frac{4}{(2\pi)^{4/3}}\,.
\end{equation}
}

\end{remark}

One of the main motivations of this paper is to fix the undetermined constant $\delta_g$ in (\ref{main}) for $g>1$. To this end we are to study the asymptotics of ${\rm det}^*\,\Delta$ when two zeros of the differential collide and the surface $\X$ degenerates to a nodal surface with two irreducible components $\X^+$ and $\X^-$.

In terminology of \cite{EMZ} we approach the principal boundary of $H_g(1, \dots, 1)$ shrinking two homologous saddle connections (i. e. geodesics, joining two colliding zeros). One can think about this situation as follows.
Let $g^+, g^-\geq 1$ be integers such that $g^++g^-=g$ and let $(\X^\pm, \omega^\pm)\in H_{g^\pm}(1, \dots, 1)$ ($2g^\pm-2$ units).
 Introduce two straight cuts, $[P_+, P_+(s)]$ and $[P_-, P_-(s)]$, of equal length $s$: one on the translation surface $\X^+$ and another on the translation surface $\X^-$ (these cuts should not    contain the conical points). Identifying each shore of the cut on the surface $\X^+$ with the corresponding  shore of the cut on the surface $\X^-$, one gets a translation surface $\X^{(s)}$ of genus $g=g^++g^-$ with $2g-2=(2g^+-2)+(2g^--2)+2$ conical points of angle $4\pi$:  $2g^+-2$ of them, $P_1^+, \dots, P_{g^+}^+$  come from the surface $\X^+$, the $2g^--2$ points, $P_1^-, \dots, P_{g^-}^+$, come from the surface $\X^-$ and the remaining two conical points, $P_r$ and $P_l$, are the end points of the cuts. One can see that the points $P_r$ and $P_l$ are joined by two homologous saddle connections of length $s$ on the surface $\X^{(s)}$, these saddle connections are just the former shores of the cuts.
 The translation surface $\X^{(s)}$ comes with holomorphic one form $\omega^{(s)}$ having simple zeros at conical points of $\X^{(s)}$ and coinciding with $\omega^\pm$ in $\X^\pm\setminus [P_\pm, P_\pm(s)]\subset \X^{(s)}$.
 So we are interest in the asymptotics of ${\rm det}^*\,\Delta^{|\omega^{(s)}|^2}$ as $s\to 0$.

 As we see the degeneration scheme we encounter here is slightly different from the
 usual one (see, e. g., \cite{Fay}, \cite{Masur}, \cite{Yamada}), where the family of degenerating Riemann surfaces is obtained from two surfaces $\X^+$ and $\X^-$ via the well known plumbing construction (one glues not the shores of the cuts as we do here but the annuli $A_\pm=\{s\leq |z_\pm|\leq 1\}\subset \X^\pm$ identifying the points $z_+$ and $z_-$ such that $z_+z_-=s$).
 Thus, one has to modify the results from \cite{Fay} (later corrected in \cite{Yamada}) concerning the asymptotical behavior of basic holomorphic objects on the degenerating Riemann surface (in particular those entering  (\ref{deftau})) in order to serve a different degeneration scheme.                    Section 2 of the present paper is devoted to this tedious but, unfortunately, indispensable task.  In this section we closely follow  Fay and Yamada, we have chosen to use a certain hybrid of their approaches in order to keep all the proofs elementary and (hopefully) a little bit more readable than their prototypes.

After this task is completed it becomes possible to calculate the asymptotics of $\tau_g(\X^{(s)}, \omega^{(s)}, \{\a_\alpha, b_\alpha\})$ from (\ref{deftau}) as $s\to 0$. The result (obtained in subsection 2.4) looks as follows
 \begin{equation}\label{astau}\tau_g(\X^{(s)}, \omega^{(s)}, \{\a_\alpha, b_\alpha\})\sim \frac{1}{\sqrt{2}}s^{1/4}\tau_{g_+}(\X^{+}, \omega^{+}, \{\a_\alpha^+, b_\alpha^+\})\tau_{g_-}(\X^{-}, \omega^{-}, \{\a_\alpha^-, b_\alpha^-\})\,;\end{equation}
here the canonical basis $\{\a_\alpha, b_\alpha\}$ on the surface $\X^{(s)}$ is  the union of the canonical basis $ \{\a_\alpha^+, b_\alpha^+\}$ on $\X^+$ and the canonical basis $\{\a_\alpha^-, b_\alpha^-\}$ on $\X^-$.

This result implies  the asymptotics
\begin{equation}\label{Main1}{\rm det}^*\,\Delta^{|\omega^{(s)}|^2}\sim \frac{\delta_g}{2\delta_{g_+}\delta_{g_-}} s^{1/2}\frac{{\rm Area}\,(\X)}{{\rm Area}\,(\X^+){\rm Area}\,(\X^-)}\,{\rm det}^*\,\Delta^+{\rm det}^*\,\Delta^-\end{equation}
and in order to fix the constant $\delta_g$ it is sufficient to get the asymptotics of  ${\rm det}^*\,\Delta^{|\omega^{(s)}|^2}$
for  some {\it special} elements $(\X^\pm_0, \omega^\pm_0)$ of $H_{g_\pm}(1, \dots, 1)$ using another method and then compare the coefficients in the two asymptotics. (It should be noted that a similar program was recently realized by R. Wentworth for the determinants of the Laplacian in the Arakelov metric in order to calculate the so-called bosonisation constants (see \cite{Wentworth})). This is done in Section 3.
The key idea (picked up by the author in a conversation with L. Hillairet) is the following: one can start (in case of even genus $g=2g_0$) with a translation surface $\X_0$ of genus $g_0$ with a cut $[P, P(s)]$ of length $s$ and glue two copies of $\X_0$ together along the cut. (So, one takes $\X^+=\X^-=\X_0$ in the above construction.)
In this symmetric situation the Laplacian $\Delta$ on the translation surface $\X^{(s)}$ is unitary equivalent to the direct sum of the two operators, $\Delta_{\cal D}$ and $\Delta_{\cal N}$, of Neumann and Dirichlet homogeneous boundary value problems in $\X_0\setminus [P, P(s)]$.
Thus, one has the relation
\begin{equation}\label{reduction}{\rm det}^*\,\Delta={\rm det}\,\Delta_{\cal D}\,{\rm det}^*\,\Delta_{\cal N}\end{equation}
(notice that the Dirichlet Laplacian has no zero modes and one does not modify its determinant here).
I turns out that the asymptotics of ${\rm det}\,\Delta_{\cal D}$ and ${\rm det}^*\,\Delta_{\cal N}$ as $s\to 0$ can be found if one makes use of a certain variant of the BFK surgery formula (see \cite{BFK}), the Wentworth lemma on the asymptotics of the Dirichlet-to-Neumann  operator  on a shrinking contour (\cite{Wentworth}) and a simple idea based on rescaling properties of the determinant of the Laplacian.
So, one can find the asymptotics ${\rm det}^* \Delta$ for a symmetric translation surface (and, therefore, for an arbitrary translation surface of genus which is an integer power of $2$); a simple trick based on BFK surgery formula reduces the general case to this symmetric one.

Beyond the scope of the present paper remains the case of another possible collision of conical points (in other words we consider here the asymptotical behavior of ${\rm det}^*\Delta$ only near a part of the principal boundary of the stratum): one can shrink a saddle connection of length $s\to 0$ which has no saddle connection homologous to it. In this case the underlying Riemann surface $\X^{(s)}$ does not degenerate (and tends to a  nonsingular Riemann surface $\X^{(0)}$; we denote by $\Delta_0$ the Laplacian on the translation surface $\X^{(0)}$) but the colliding zeros form a single zero of multiplicity two (a conical point of the angle $6\pi$). It is relatively easy to show that in this case the asymptotics of $\tau_g$ has the form
$$ \tau_g(\X^{(s)}, \omega^{(s)}, \{\a_\alpha, b_\alpha\})\sim s^{1/36}\tilde \tau_{g}(\X^{(0)}, \omega_0, \{\a_\alpha, b_\alpha\}),$$
where $\tilde \tau_g$ is an analog of the function $\tau_g$ for the stratum $H_g(2, 1, \dots, 1)$ ($2g-4$ units; see \cite{DG} for definitions). This (together with results from \cite{DG}) leads to the asymptotics ${\rm det}^*\Delta\sim C_g s^{1/18}{\rm det}^*\Delta_0$ with the unknown constant $C_g$. Finding this constant presents an interesting open problem. Even more complicated looks the problem of finding the asymptotics of ${\rm det}^* \Delta$ at the boundary of a general stratum $H_g(k_1, \dots, k_M)$, at the moment we see no reasonable approach to it.

Finally we notice that similar problems for hyperbolic metric of constant curvature were studied by S. Wolpert (\cite{Wolp}) and R. Lundelius (\cite{Lund}), the case of Arakelov metric (with curvature given by the Bergman 2-form) was considered in \cite{Jorg} and \cite{Went1, Went2} (the complete results were recently obtained by R. Wentworth in \cite{Wentworth}).
We think that the case of a metric with curvature  concentrated at a finite set  (considered in the present work) forms a natural complement to these results filling the right hand side of the picture (if one puts the constant curvature metric at the left hand side and the Arakelov metric in the center).

The author is grateful to R. Wentworth for explaining some subtle details from \cite{Wentworth} and clarification of the reason of divergence between the results of \cite{Fay} and \cite{Yamada}; the author also thanks  L. Hillairet for generous sharing of his ideas on spectral theory of translation surfaces and D. Korotkin for numerous useful discussions.

\section{Families of degenerating surfaces and asymptotical formulas}

 We construct several one-parametric families of Riemann surfaces degenerating as the parameter tends to zero.

 Let $\X^+$ and $\X^-$ be two compact Riemann surfaces of genus $g^+$
 and
$g^-$, $g^\pm\geq 0$. Choose points $P_{\pm}\in \L^{\pm}$ and their
 open
neighborhoods $D^\pm\subset \L^{\pm}$ such that for a certain choice of
holomorphic local parameters $\l^\pm$ on $\L^\pm$ one has
$D^{\pm}=\{P\in \L^\pm :|\l^\pm(P)|<1\}$ and  $\l^{\pm}(P^\pm)=0$. Define the map
 $\l:D^+\cup
D^-\to {\mathbb C}$ setting $\l(P)=\l^\pm(P)$ if $P\in D^\pm$.

Using these data we construct three families of degenerating Riemann surfaces
of genus $g^-+g^+$.

{\bf Case I.} Let $s$ be a complex number, $|s|<1$ and let $P^\pm(s)$ be the
points
 in
$D^{\pm}$ such that $\l(P^\pm(s))=s$.

Cut the discs $D^{\pm}$ along the (oriented) straight segments $[P^\pm,
P^\pm(s)]$ and glue the surfaces $\L^+$ and $\L^-$ along these cuts
identifying a point $P$ on the left shore of the "$+$"-cut with the
 point $Q$
($\l^+(P)=\l^{-}(Q)$) on right shore of the "$-$"-cut and vice versa;
 the
resulting topological real 2-d surface can be turned into a compact
 Riemann
surface $\L_s$ of genus $g=g^-+g^+$ in a usual way (one chooses the
 local
parameter near the left endpoint P of the cut as
 $\zeta(Q)=\sqrt{\l(Q)}$,
near the right endpoint $P(s)$ the local parameter is
$\zeta(Q)=\sqrt{\l(Q)-\l(P(s))}$, the choice of the local parameter at other points of $\L_s$ is obvious).

{\bf Case Ia.} This family is constructed similarly to Cases I, the only
difference is the position of cuts inside the disks $D^\pm$: choose
 a complex number $t$, $|t|<1$ and introduce the cuts inside the discs $D^\pm$ connecting the points $\l=\sqrt{t}$ and $\l=-\sqrt{t}$; after
  the same gluing  of the shores of these cuts as in case I we get the family  $\L_t$  of degenerating  compact Riemann surfaces.

{\bf Case Ib.} This family is obtained similarly to Cases I and Ia, but
instead of gluing the disks along the cuts we use the standard "plumbing
construction" (see \cite{Fay}). Choose $t$, $|t|<1$ delete from the discs
$D^\pm$ the smaller discs $|\l^\pm|\leq |t|$ and glue the obtained annuli,
$A^\pm$, identifying points $P\in A^+$ and $Q\in A^-$ such that
$\l^+(P)\l^-(Q)=t$. After this gluing the surfaces $\L^\pm$ turn into a
single Riemann surface $\L_t'$ of genus $g^-+g^+$.

In what follows we derive asymptotical formulas (as $s\to 0$) for basic
holomorphic objects (the normalized holomorphic differentials, the canonical
meromorphic differential, the prime-form, etc) on the Riemann surfaces
constructed in case I.

The asymptotical formulas (as $t\to 0$) for case Ib  were first derived in
\cite{Fay}.  In \cite{Yamada} it was claimed that all the formulas from
\cite{Fay} are incorrect and new ones were proved. Our analysis (in
particular, see  Example 1 below) shows that formulas from \cite{Fay} (as
well as Fay's proofs of these formulas) are applicable in case Ia. As it was
explained to us by Richard Wentworth (private communication) Fay in fact
makes a mistake when considering case Ib: his  "pinching parameter" depends
in its turn on deformation parameter and this results in additional terms in
asymptotical expansions which were lost in \cite{Fay}. In case Ia  the
pinching parameter is independent of deformation parameter and Fay's scheme
works perfectly.

The case of our concern, I, is very similar to case Ia  (the pinching
parameter,  $\l$ in equation (\ref{MAIN}) below, is independent of the
deformation parameter $s$) and we give here the proofs of all the
asymptotical formulas for it. Mainly we use the methods similar to those of
Fay (where they are applicable); although we have chosen to follow the pretty
elementary analytical methods of Yamada (avoiding Grauert's theorem and sheaf
cohomologies from \cite{Fay}, \cite{Masur}) when introducing a holomorphic
family of Abelian differentials on $\L_s$ and studying the analytical
properties of the coefficients in the Laurant expansions in the pinching
zone.

\subsection{Two examples in genus 0}

{\bf Canonical meromorphic bidifferential $W$.}
Recall that to any compact Riemann surface $\X$ of genus $g$ with a chosen canonical basis of cycles $\{a_\a, b_\a\}$ on it one associates the so-called canonical meromorphic bidifferential $W(\,\cdot\,, \,\cdot\,)$, which
\begin{itemize}
\item is a meromorphic one-form with respect to each argument,
\item is symmetric, i. e. $W(P, Q)=W(Q, P)$
\item has a single pole at the diagonal $P=Q$ and
$$W(z(P), z(Q))\sim \frac{dz(P)dz(Q)}{(z(P)-z(Q))^2}+\frac{1}{6}S_B(z(Q))dz(P)dz(Q)+o(1)$$
as $P\to Q$, where $S_B$ is the Bergman projective connection (see \cite{Fay}).
\item satisfies  $\oint_{a_\a}W(P,\, \cdot\,)=0$ for any $\P\in \X$ and $\a=1, \dots, g$.
\end{itemize}
(In case $g=0$ the last condition is void.) The canonical meromorphic bidifferential is related to the prime form via the equation
 $W(P, Q)=d_Pd_Q\log E(P, Q)$ (see \cite{Fay}).

{\bf Example 1: case Ia.}
We start with the following simple statement. Let $\L$ be the two-fold
branched covering of the Riemann sphere ${\mathbb P}^1$ with branch
 points
$\l_1$, $\l_2$. Let $P\in \L$ and $\l$ be the projection of $P$ on
 ${\mathbb
P}^1$. Then the map
$$P\mapsto \delta=\sqrt{\frac{\l-\l_1}{\l-\l_2}}$$
is the biholomorphic isomorphism of $\L$ and ${\mathbb P}^1$. Applying
 to
$\delta$ the fractional linear transformation $\delta\mapsto \gamma=
\frac{\l_2-\l_1}{\delta-1}+\l_2$, we get the isomorphism

\begin{equation}\label{1}P\mapsto
 \gamma=\l+\sqrt{(\l-\l_1)(\l-\l_2)}\end{equation}
 of $\L$ and   ${\mathbb
P}^1$ which is more convenient for our future purposes.

Now let $t>0$ and $\l_1=-\sqrt{t}$, $\l_2=\sqrt{t}$. When $t\to 0$ the
Riemann sphere $\L_t$ degenerates to the singular Riemann surface
 with two
components, Riemann spheres $S^+$ and $S^-$, attached to each other at
 the
point $0$. So, our situation is exactly the one described in Case Ia.

Let $W_t(\cdot, \cdot)$, $W_-(\cdot, \cdot)$ and
 $W_+(\cdot,
\cdot)$ be the canonical meromorphic bidifferentials on $\L_t$, $S^-$
 and
$S^+$ respectively.

Then the following asymptotics holds:

\begin{equation}\label{Fay1}W_t(\l(P), \l(Q))=
\begin{cases} W_\pm(\l(P), \l(Q))+\frac{t}{4}W_\pm(\l(P),
 0^\pm)W_\pm(\l(Q), 0^\pm)+O(t^2), \ \ \ \  \text{if}\ \ P, Q\in
 S^{\pm}\\
-\frac{t}{4}W_\pm(\l(P), 0^\pm)W_\mp(\l(Q), 0^\mp)+O(t^2) \ \
 \ \ \
\text{if}\ \ \ \ P\in S^\pm, Q\in S^\mp
\end{cases}\end{equation}

(This asymptotics  (with the minus sign in the last line lost) was mistakenly stated in (\cite{Fay}, formula (49), p. 41) for two
Riemann surfaces glued via plumbing construction (Case Ib), however, being false in Case Ib, it is true in Case Ia.)

Let $P, Q$ be two points of the covering $\L_t$ lying on the same sheet
 (say
$S_+$) with projections $\l$ and $\zeta$; assume for simplicity that $\l$
 and
$\zeta$ are real and positive.

Using the uniformization map (\ref{1}), one can write the following
asymptotics for the canonical meromorphic differential on $\L_t$:
$$W_t(\l,
\zeta)=\frac{d\gamma(\l)d\gamma(\zeta)}{(\gamma(\l)-\gamma(\zeta))^2}=\frac{(1+\frac{\l}{\sqrt{\l^2-t}})(1+\frac{\zeta}{\sqrt{\zeta^2-t}})}
{[\l-\zeta+\sqrt{\l^2-t}-\sqrt{\zeta^2-t}]^2}d\l\,d\zeta=$$
\begin{equation}\label{2}=\frac{d\l\,d\zeta}{(\l-\zeta)^2}+\frac{t}{4\l^2\zeta^2}d\l\,d\zeta+O(t^2)d\l\,d\zeta\end{equation}
as $t\to 0+$ which agrees with Fay's formula (49).

(We remind the reader that the canonical bidifferential $W_+$ on
 $S_+$
(as well as on $S_-$) is $\frac{d\l\,d\zeta}{(\l-\zeta)^2}$ and, therefore,
$W_\pm(\l(P), 0^\pm)W_\pm(\l(Q),
 0^\pm)=\frac{1}{\l^2\zeta^2}d\l
d\zeta$.)

If $P\in S^+$, $Q\in S^-$ then all the "$\zeta$"- square roots in
 (\ref{2})
change their sign and we  arrive at the second case of Fay's expansion (\ref{Fay1}).

{\bf Example 2: Case Ib.} This is a rather elementary simplification of Yamada's Example 1 (\cite{Yamada}, pp. 140-142), the author thanks D. Korotkin for pointing it out to him.

Let $S^+=S_v$ and $S^-=S_w$ be two Riemann spheres with standard coordinates $v$ and $w$ in $S_{v, w}\setminus \{\infty\}$. Let also $\zeta=1/w$ be the local parameter near the point at infinity of $S_w$.
Glue $S_v\setminus \{|v|<t\}$ and $S_w\setminus\{|\zeta|<t\}$ together identifying the points $v\in \{t\leq |v|\leq 1\}\subset S_v$ and
$\zeta\in  \{t\leq|\zeta|\leq 1\}\subset S_w$ such that $v\zeta=v/w=t$. We get a Riemann surface $\X_t$ of genus $0$.
It is easy to write the uniformization map $\X_t\to S_z$, where $S_z$ is the Riemann sphere with standard coordinate $z$ in $S_z\setminus\{\infty\}$.

Namely, define the map $z: S_v\setminus \{|v|<t\}\cup S_w\setminus\{|\zeta|<t\}\to S_z$ via
$z(v)=v$ for $v\in S_v\setminus \{|v|<t\}$ and $z(w)=t w$ for $w\in S_w\setminus\{|\zeta|<t\}$. Obviously, the relation $v/w=t$ implies $z(v)=z(w)$, therefore, the map $z$ gives rise to a biholomorphic map $\X_t\to S_z$.

One has the following obvious relations for the canonical meromorphic bidifferentials $W_t$, $W_\pm$ on $\X_t$ and $\X_\pm$.
\begin{equation}\label{eq111}
W_t(z_1, z_2)=\frac{dz_1dz_2}{(z_1-z_2)^2}=\frac{dv_1dv_2}{(v_1-v_2)^2}=W_+(v_1, v_2),
\end{equation}
if $v_1, v_2\in S^+\setminus\{|v|\leq 1\}$;

\begin{equation}\label{eq222}
W_t(z_1, z_2)=\frac{dz_1dz_2}{(z_1-z_2)^2}=\frac{d(tw_1)d(tw_2)}{(tw_1-tw_2)^2}=\frac{dw_1dw_2}{(w_1-w_2)^2}=W_-(w_1, w_2),
\end{equation}
if $w_1, w_2\in S^-\setminus\{|\zeta|\leq 1\}$;

\begin{equation}\label{eq333}
W_t(z_1, z_2)=\frac{dz_1dz_2}{(z_1-z_2)^2}=\frac{dv d(tw)}{(v-tw)^2}=t\frac{dvdw}{v^2}+O(t^2)
\end{equation}
as $t\to 0$, if $v\in S^+\setminus\{|w|\leq 1\}$ and $w\in  S^-\setminus\{|\zeta|\leq 1\}$ in complete agreement with Yamada's asymptotical formulas
for the case Ib:
\begin{equation}\label{Yama1}
W_t(z_1, z_2)=\begin{cases}
W_+(v_1, v_2)+t^2S_B(\zeta)|_{\zeta=0}W_+(v_1, 0)W_+(v_2, 0)+O(t^3)\ \ \ {\text for}\ \ v_1, v_2\in S^+\setminus\{|v|\leq 1\};\\
W_-(w_1, w_2)+t^2S_B(v)|_{v=0}W_-(w_1, \infty)W_-(w_2, \infty)+O(t^3)\ \ \ {\text for}\ \ w_1, w_2\in S^-\setminus\{|\zeta|\leq 1\};\\
-tW_+(v, 0)W_-(w, \infty)+O(t^2) \ \ \ {\text for}\ \ v\in S^+\setminus\{|v|\leq 1\}, w\in S^-\setminus\{|\zeta|\leq 1\}
\end{cases}
\end{equation}
(see \cite{Yamada}, formula (15) on p. 122; it should be noted that for coordinates $v$ and $\zeta$ on the Riemann sphere one has $S_B(v)=S_B(\zeta)\equiv 0$).

\subsection{Asymptotical formulas.}

Here we deal with Case I, assuming that the genera of the surfaces $\L^\pm$
are greater than zero.

Denote the part of the Riemann surface $\L_s$ which came from the discs $D^\pm$ after the gluing procedure by ${\cal U}$. The domain ${\cal U}$ is an open (topological) annulus
and the map $\l$ can be considered as defined on ${\cal U}$. The map
\begin{equation}\label{cover}\l: {\cal U}\rightarrow \{\l\in {\mathbb C} : |\l|<1\}\end{equation}
defines a two-sheeted covering of the disc $\{|\l|<1\}$ ramified over $\l=0$ and $\l=s$, whereas the map
\begin{equation}\label{unifor}{\cal U}\ni P\mapsto X=\l-\frac{s}{2}+\sqrt{\l(\l-s)}\end{equation}
is a well-defined biholomorphic bijection (of course, the value of the square
root depends on to which disk, $D^+$ or $D^-$, the point $P$ belongs; one
also has to fix a branch of the square root, say, for the disk $D^+$ with the
cut between $0$ and $s$, there are two choices and we make one once and
forever).

(It should be noted that  map (\ref{unifor})  (being appropriately extended)
 uniformizes the two-sheeted covering of the Riemann sphere branched over the points $0$ and $s$. The image of the point at
infinity of the first sheet is $\infty$, the image of the point at infinity of the second one is $0$.)

For sufficiently small $s$ the annulus
$${\mathbb A}_s=\{P : \frac{|s|^2}{4}< |X|<1\}$$
belongs to ${\cal U}$. Moreover, the boundary curve $|X|=1$ lies in a small
vicinity of the circle $|\l|=1/2$ of the "$+$"-sheet of the covering
(\ref{cover}), whereas the boundary curve $|X|=|s|^2/4$ lies in a small
vicinity of the circle $|\l|=1/2$ of the "$-$"-sheet.

The following two lemma are analogs of Yamada's Theorem 1 and Lemma 1
(\cite{Yamada}, p. 116) for the family $\L_s$. We follow the proofs of Yamada
making necessary (in fact, rather minor)  modifications.

\begin{lemma}\label{Yama}
Let $v_\pm$ be holomorphic differentials on $\L^\pm$. There exists a
holomorphic differential $w_s$ on $\L_s$ such that for any $\rho$,
$\sqrt{|s|}<\rho<1$   holds the inequality
\begin{equation}
||w_s-v_+||_{\Omega^+_\rho}+||w_s-v_-||_{\Omega^-_\rho}\leq C(\rho)|s|,
\end{equation}
where
$$\Omega^+_\rho=\L^+\setminus \{P\in D^+ : |X(P)|\leq \rho\}\,$$
$$\Omega^-_\rho=\L^-\setminus  \{P\in D^- : |X(P)|\geq |s|^2/(4\rho)\}\,.$$
Here as usual, the $L_2$-norm of a one-form in a subdomain $\Omega$ of a Riemann surface is defined via $$||u||_\Omega=\int\int_\Omega u\wedge \overline{*u}\,.$$
\end{lemma}

{\bf Remark.} The curves $|X|=\sqrt{|s|}$ and $|X|=|s|^{3/2}/4$ belong to
small vicinities of the circles $|\l|=\sqrt{|s|}/2$ lying on the "$+$" and
"-" sheets of the covering (\ref{cover}) respectively.

{\bf Proof.}

Let $\int_0^\l u_+=\sum_{n=1}^\infty \alpha_n\l^n$ near $P^+$; after passing
to coordinate $X$,
$$\l=\frac{X}{2}+\frac{s}{2}+\frac{s^2}{8X},$$
we get
$$f_+(\l)=\int_0^\l u^+=\sum_{n=1}^\infty a_n^+(s)X^n+a_0(s)+\sum_{n=-\infty}^{-1}a_n^{-}(s)X^n,$$
where
$$a_n^+(s)=\alpha_n(1/2^n+O(s)); \ \ a_0(s)=O(s); \ \
a^-_n(s)=O(s^{2|n|})\,,$$ as $s\to 0$.

Analogously, from the expansion the expansion $f_-(\l)=\int_0^\l
u_-=\sum_{n=1}^\infty \beta_n\l^n$ near $P^-$ one gets
$$\int_0^\l u_-=\sum_{n=1}^\infty b_n^+(s)X^n+b_0(s)+\sum_{n=-\infty}^{-1}b_n^{-}(s)X^n,$$
where
$$b_n^+(s)=\beta_n(1/2^n+O(s)); \ \ b_0(s)=O(s); \ \
b^-_n(s)=O(s^{2|n|})\,,$$ as $s\to 0$.

Now,  \cite{Yamada}, we are to construct a sequence, $\Phi^{(k)}_s$, of
$C^1$-forms on $\L_s$ coinciding with $v_\pm$ in $\Omega_\rho^\pm$ and such
that
\begin{equation}\label{in2}||\Phi^{(k)}_s-i*\Phi^{(k)}_s||^2\leq O(s^2)+1/k\,
.\end{equation}

For harmonic function $h_s$ in the annulus $\{|s|^2/4\rho\leq |X|\leq \rho\}$
with boundary values $f_-$ and $f_+$ one has the relation
\begin{equation}\label{in0}\frac{1}{2\pi}\int\int_{|s|^2/4\rho\leq|X|\leq
\rho}(|\partial_Xh_s|^2+|\partial_{\bar X}h_s|^2)\frac{|dX\wedge \overline
{dX}|}{2}=$$$$=\sum_{n=1}^\infty\frac{n|b^-_{-n}-a^-_{-n}|^2}{\rho^{2n}-(\frac{|s|^2}{4\rho})^{2n}}+\sum_{n=1}^\infty\frac
{n|b^+_n-a^+_n|^2}{\rho^{-2n}-(\frac{|s|^2}{4\rho})^{-2n}}+\frac{|b_0-a_0|^2}{2\log(\frac{\rho^2}{|s|^2/4})}=O(s^2)\,.
\end{equation}

It can be shown (say, via polynomial interpolation along radii directions)
that one can change the function $h_s$ in small vicinities of boundary
circles $|X|=\rho$ and $|X|=|s^2|/4\rho$ obtaining the function $h_s^{(k)}$
such that

\begin{equation}\label{in1}\int\int_{|s|^2/4\rho\leq|X|\leq
\rho}(|\partial_X(h_s-h_s^{(k)})|^2+|\partial_{\bar
X}(h_s-h_s^{(k)})|^2)\frac{|dX\wedge \overline
{dX}|}{2}\leq\frac{1}{k}\,\end{equation} and the $1$-form
\begin{equation}
\Phi_s^{(k)}=\begin{cases} v_{\pm} \ \ \text{in}\ \ \Omega_\rho^\pm, \\
d(h_s^{(k)})\ \ \text{in}\ \ \L_s\setminus (\Omega_\rho^+\cup\Omega_\rho^-)
\end{cases}
\end{equation}
is $C^1$-smooth. Since the operator $\text{Id}-i*\ $ kills the $(1,0)$-forms,
the inequality (\ref{in2}) follows from (\ref{in0}) and (\ref{in1}).

Decomposing $(\text{Id}-i*)\Phi_s^{(k)}$ into ($L_2$-orthogonal!) sum of a
harmonic one-form $\omega_h$, an exact form $\omega_e$ and a co-exact form
$\omega^{*}_e$ (see \cite{Ahlfors}, Chapter V; here "exact form" means a form
belonging to the $L_2$-closure of the space of smooth exact forms), we
observe that the left part of the equation
$$\Phi_s^{(k)}-\omega_e=i*\Phi_s^{(k)}+\omega_h+\omega_e^*$$
is a closed form, whereas its left part is co-closed, therefore, both are
harmonic by virtue of Weyl's Lemma (see \cite{Ahlfors}, Chapter V).

Now, applying to the harmonic form $\Phi_s^{(k)}-\omega_e$ the operator
$\frac{1}{2}(\text{Id}+i*)$ one gets a holomorphic one-form
$$\Psi^{(k)}_s=\frac{1}{2}(\text{Id}+i*)[\Phi_s^{(k)}-\omega_e]\,.,$$
which coincides with $v_\pm+\frac{1}{2}(\text{Id}+i*)\omega_e$ in
$\Omega^\pm_\rho$. Therefore,
\begin{equation}\label{in4}
||\Psi^{(k)}_s-v_+||^2_{\Omega^+_\rho}+||\Psi^{(k)}_s-v_-||^2_{\Omega^-_\rho}\leq\frac{1}{4}||\omega_e+i*\omega_e||\leq
\frac{1}{2}||\omega_e||\leq\frac{1}{2}||\Phi_s^{(k)}-i*\Phi_s^{(k)}||
\end{equation}
and
$$||\Psi^{(k)}_s-v_+||^2_{\Omega^+_\rho}+||\Psi^{(k)}_s-v_-||^2_{\Omega^-_\rho}\leq
O(s^2)+\frac{1}{k}$$ by virtue of (\ref{in2}).

Choosing from the sequence $\{\Psi_s^{(k)}\}_{k\geq 1}$ a converging
subsequence (uniform $L_2$-boundedness of holomorphic forms on a compact
Riemann surface implies uniform boundedness of their coefficients) and
passing to the limit $k\to \infty$ we get a holomorphic 1-form $w_s$ with all
the needed properties. $\square$
\begin{remark}\label{Polusa}{\rm
Actually a stronger variant of Lemma 1 is true: the differentials $v_\pm$ can
be meromorphic with poles lying outside of $D^\pm$. In this case the
differential $w_s$ is also meromorphic and have the same singularities as
$v_\pm$.}
\end{remark}

Now choose on $\L^\pm$ a canonical basis of cycles $\{a_\alpha^\pm,
b_\alpha^\pm\}_{\alpha=1, \dots, g^\pm}$ such that none of the cycles
intersects the disk $D^\pm$. Let also $\{u_\alpha^\pm\}_{\alpha=1, \dots,
g^\pm}$ be the corresponding basis of normalized differentials.

The set of cycles $\{a_\alpha, \b_\alpha\}_{\alpha=1, \dots,
g^++g^-}=\{a^+_1, \dots, a^+_{g^+},a^-_1, \dots,a^-_{g^-}; b^+_1, \dots,
b^+_{g^+}, b^-_1, \dots, b^-_{g^-}\}$ forms a canonical basis on the Riemann
surface $\L_s$.
 Let $\{ v^{(s)}_\alpha \}_{\alpha=1, \dots, g^-+g^+}$
be the corresponding basis of normalized holomorphic differentials on $\L_s$.

Let also $w^{(s)}_\alpha$ be a holomorphic one form on $\L_s$ which is
constructed in Lemma 1 when one takes $(v_+, v_-)=(v_\alpha^+, 0)$  for
$\alpha=1, \dots, g^+$ and $(v_+, v_-)=(0, v^-_{\alpha-g^+})$ for
$\alpha=g^++1, \dots, g^++g^-$.

The corresponding $a$-period matrix ${\mathbb
P}=||\oint_{a_\alpha}w^{(s)}_\beta||_{\alpha, \beta=1, \dots, g^++g^-}$
satisfies
$${\mathbb P}=I_{g^++g^-}+O(s)$$
as $s\to 0$ due to Lemma 1. This immediately implies the following lemma.

\begin{lemma} The basis $\{v^{(s)}_\alpha\}_{\alpha=1, \dots, g^++g^-}$ of
normalized holomorphic differentials on $\L_s$ satisfies
\begin{equation}
(v^{(s)}_1, \dots, v^{(s)}_{g^++g^-})=(I_{g^-+g^+}+O(s))(w^{(s)}_1, \dots,
w^{(s)}_{g^++g^-}),
\end{equation}
in particular, all the differentials $v^{(s)}_\alpha$ are uniformly (with
respect to $s$) bounded in, say,  $\L_s\setminus \{P\in \L_s, |\l(P)|<1/4\}$.
\end{lemma}

{\bf Laurent expansion for basic holomorphic differentials.} Writing the
differential $v^{(s)}_\alpha$ as $v^{(s)}_\alpha(X)dX$ in the local parameter
$X=\l-\frac{s}{2}+\sqrt{\l(\l-s)}$ and expanding the coefficient
$v^{(s)}_\alpha(\cdot)$ in the Laurent series in the annulus $|s|^2/4<|X|<1$,
one gets
\begin{equation}\label{Lor}v_\alpha^{(s)}(X)dX=(\sum_{n>o}\gamma_{-n}(s)X^{-n}+\sum_{n\geq
0}\gamma_n(s)X^n)dX\,.\end{equation} Observe that
$dX=\frac{Xd\l}{\sqrt{\l(\l-s)}}$ and for $n\geq0$ one has
\begin{equation}\label{pol}X^ndX=\frac{\left(\l-s/2+\sqrt{\l(\l-s)}\right)^{n+1}}{\sqrt{\l(\l-s)}}d\l=\left\{\sum_{k=0}^{n+1}p_k(s)\l^k+
\frac{1}{\sqrt{(\l(\l-s)}}\sum_{k=0}^{n+1}q_k(s)\l^k
\right\}d\l\end{equation} with some polynomials $p_k(s), q_k(s)$. On the
other hand, since
$$(\l-s/2+\sqrt{\l(\l-s)})(\l-s/2-\sqrt{\l(\l-s)})=s^2/4\,,$$
 for $n>0$ one has
$$X^{-n}dX=\frac{4^n}{s^{2n}}\frac{\left(\l-s/2-\sqrt{\l(\l-s)}\right)^n\left(\l-s/2+\sqrt{\l(\l-s)}\right)}{\sqrt{\l(\l-s)}}d\l=$$
\begin{equation}\label{otr}=\frac{1}{s^{2n-2}}\left\{\sum_{k=0}^{n-1}\tilde{p}_k(s)\l^k+
\frac{1}{\sqrt{\l(\l-s)}}\sum_{k=0}^{n-1}\tilde{q}_k(s)\l^k \right\}d\l
\end{equation}
with some polynomials $\tilde{p}_k(s), \tilde{q}_k(s)$.

For $n>0$ one has
$$\gamma_{-n}(s)=\frac{1}{2\pi i}\int_{|X|=|s|^2/4}
v_\alpha^{(s)}(X)X^{n-1}dX=\frac{1}{2\pi i}\int_{\Gamma_-}
v_\alpha^{(s)}(\l)\left(\l-s/2+\sqrt{\l(\l-s)}\right)^{n-1}\,d\l=$$
\begin{equation}\label{otr1}=\int_{\Gamma_-} O(1)\times O(s^{2n-2})d\l= O(s^{2n-2})
\end{equation}
as $s\to 0$ (the contour $\Gamma_-$ over which goes the last integration lies
in a small vicinity of the circle $|\l|=1/2$ of the "-"-sheet; the factor
$v_\alpha^{(s)}(\l)$ is uniformly bounded on this contour with respect to $s$
by virtue of Lemma 2).

In the same manner for $n\geq0$ one has
\begin{equation}\label{pol1}\gamma_n(s)=\frac{1}{2\pi
i}\int_{|X|=1}\frac{v_\alpha^{(s)}(X)}{X^{n+1}}dX=\frac{1}{2\pi
i}\int_{\Gamma_+}\frac{v_\alpha^{(s)}(\l)d\l}{\left(\l-s/2+\sqrt{\l(\l-s)}\right)^{n+1}}=O(1)\end{equation}
(The contour $\Gamma_+$ lies in a small vicinity of the circle $|\l|=1/2$ of
the $+$-sheet, the factor $v^{(s)}_\alpha(\l)$ is uniformly bounded by virtue
of Lemma 2, the denominator of the integrand is close to $1$.)

Now from (\ref{Lor}), (\ref{pol}) and (\ref{otr}) together with the estimates
(\ref{pol1}) and (\ref{otr1}) one gets the expansion
\begin{equation}\label{MAIN}
v_\alpha^{(s)}(\l)d\l=\sum_{k=0}^\infty
a_k(s)\l^k\,d\l+\frac{1}{\sqrt{\l(\l-s)}}\sum_{k=o}^\infty b_k(s)\l^k\,d\l,
\end{equation}
where the coefficients $a_k$, $b_k$ are {\it analytic near} $s=0$. This
expansion is valid in the  zone $\{|s|^2/4<|X|<1\}$ (the latter for small $s$
is close to the set $\{P\in \L_s : |\l(P)|\leq 1/2\}$).

\begin{remark} {\rm Expansion (\ref{MAIN}) is a complete analog of Fay's expansion
stated on page 40 of  \cite{Fay} for deformation family Ib. However, it is
important here that in (\ref{MAIN}) the parameter $\l$ is
$s$-independent whereas in expansion from \cite{Fay} the pinching parameter
$\chi$ depends on deformation parameter. The latter fact was missed by Fay
when he wrote his asymptotical expansions (in particular, his last formula on
page 40 of \cite{Fay} should contain more terms at the right hand side)
(\cite{Went}).}
\end{remark}

{\bf Main asymptotical formulas for basic holomorphic differentials and the
canonical meromorphic bidifferential.} Let $W, W_\pm$ be the canonical
meromorphic bidifferentials on  $\L_s$ and $\L^\pm$ respectively.

\begin{theorem}\label{Abdiff}
For $\alpha=1, \dots, g^+$ one has the asymptotics as $s\to 0$
\begin{equation}\label{adif1}
v^{(s)}_\alpha(P)=
\begin{cases}
u^+_\alpha(P)+\frac{s^2}{16}u^+_\alpha(P_+)W_+(P, P_+)+o(s^2)\ \ \ \text{if}
\ \ P\in L^+\setminus D^+\subset \L_s\\
-\frac{s^2}{16}u^+_\alpha(P_+)W_-(P, P_-)+o(s^2)\ \ \ \ \text{if}\ \ P\in
\L^-\setminus D^-\subset \L_s\,.
\end{cases}
\end{equation}
For $\alpha=g^++k$, $k=1, \dots, g^-$ one has
\begin{equation}\label{adif2}
v^{(s)}_\alpha(P)=
\begin{cases}
u^-_k(P)+\frac{s^2}{16}u^-_k(P_-)W_-(P, P_-)+o(s^2)\ \ \ \text{if}
\ \ P\in L^-\setminus D^-\subset \L_s\\
-\frac{s^2}{16}u^-_k(P_-)W_+(P, P_+)+o(s^2)\ \ \ \ \text{if}\ \ P\in
\L^+\setminus D^+\subset \L_s\,.
\end{cases}
\end{equation}

Here the values of differentials at the points $\P_\pm$ are calculated in the
local parameter $\l$, the values of differentials at $P\in\L^\pm\setminus
D^\pm\subset \L_s$ are calculated in an arbitrary local parameter inherited
from $\L^\pm$ (of course, the same for the l. h. s. and the r. h. s.)
\end{theorem}

\begin{theorem}\label{bergas}
For the canonical meromorphic differential on $\L_s$ one has the following
asymptotics as $s\to 0$:
\begin{equation}\label{ba}
W(R, S)=\begin{cases} W_+(R, S)+\frac{s^2}{16}W_+(R, P_+)W_+(S, P_+)\ \
\text{if} \ \ R, S\in\L^+\setminus D^+\subset \L_s,\\
-\frac{s^2}{16}W_+(R, P_+)W_-(S, P_-)\ \ \ \text{if}\ \ R\in \L^+\setminus
D^+\subset \L_s; \ S\in \L^-\setminus D^-\subset \L_s, \\
W_-(R, S)+\frac{s^2}{16}W_-(R, P_-)W_-(S, P_-)\ \ \text{if}\ \ R, S\in
\L^-\setminus D^-\subset\L_s\, .
\end{cases}
\end{equation}
\end{theorem}
{\bf Proof.} Observe that $\lim_{s\to 0}\sqrt{\l(P)(\l(P)-s)}=\pm \l(P)$ if
$\P\in D^\pm\setminus [0, s]\subset \L_s$. Let $\alpha=1, \dots, g^+$. Taking
two points in ${\cal U}$ with $\l(P)=\l$  and sending $s\to 0$ in
(\ref{MAIN}), one gets
$$u_\alpha^+(\l)d\l=\left(\sum_{k=0}^\infty a_k(0)\l^k+\sum_{k=0}^\infty
b_k(0)\l^{k-1}\right)\,d\l$$ for the point on the "$+$"-sheet
 and
$$0=\sum_{k=0}^\infty a_k(0)\l^k-\sum_{k=o}^\infty b_k(0)\l^{k-1}$$ for the
point on the "$-$"-sheet. This implies the relations
\begin{equation}\label{rel1}
b_0(0)=0
\end{equation}
and
\begin{equation}\label{rel2}
\frac{u^+_\alpha(P_+)}{2}=a_0(0)=b_1(0).
\end{equation}
For $P\in D^+$ one has
$$\frac{1}{s}(v_\alpha^{(s)}-v_\alpha^{(0)})=\sum_{k\geq
0}^\infty\frac{a_k(s)-a_k(0)}{s}\l^k\,d\l+$$
\begin{equation}\label{ee1}
=\sum_{k\geq
0}\left\{\frac{b_k(s)-b_k(0)}{s}\frac{\l^k}{\sqrt{\l(\l-s)}}+b_k(0)\l^{k-1}\frac{\frac{\l}{\sqrt{\l(\l-s)}}-1
}{s} \right\}\,d\l=\end{equation}
$$=\left\{\sum_{k=0}^\infty a_k'(0)\l^k+\sum_{k=0}^\infty
b_k'(0)\l^{k-1}+\frac{1}{2}\sum_{k=0}^\infty b_k(0)\l^{k-2}+O(s)
\right\}\,d\l\,.$$ Since $b_0(0)=0$, the limit of the left hand side of
(\ref{ee1}) as $s\to 0$ is a meromorphic differential on $\L^+$ with a single
pole at $P_+$, therefore, it is a holomorphic differential, i. e.
\begin{equation}\label{rel3}
b_0'(0)+\frac{1}{2}b_1(0)=0\,.
\end{equation}
Moreover, since all the $a$-periods of this differential vanish it equals to
zero.

Then, again for a point on the "$+$"-sheet, we have
$$\frac{1}{s^2}(v_\alpha^{(s)}-v_\alpha^{(0)})=\frac{1}{s^2}\left[\sum_{k\geq
0}(a_k(0)+sa_k'(0)+\frac{s^2}{2}a_k''(0)+O(s^3))\l^k+\right.$$
$$\left.\sum_{k\geq
0}(b_k(0)+sb_k'(0)+\frac{s^2}{2}b_k''(0)+O(s^3))\l^{k-1}(1+\frac{s}{2\l}+\frac{3}{8}\frac{s^2}{\l^2}+O(s^3))-\sum_{k\geq
0}a_k(0)\l^k-\sum_{k\geq 0}b_k(0)\l^{k-1}\right]\,d\l$$

 Since $s$-linear term
in the braces vanishes, the limit of this expression as $s\to 0$ equals to
$$\left[\sum_{k=0}^\infty
\frac{a_k''(0)}{2}\l^k+\frac{b_k''(0)}{2}\l^{k-1}+\frac{3}{8}b_k(0)\l^{k-3}+\frac{b_k'(0)}{2}\l^{k-2}\right]\,d\l.$$
Thus the limit is a meromorphic differential on $\L^+$ with a single pole of
{\it the second order} ($b_0(0)=0$!); the corresponding Laurent coefficient
is
$$\frac{3}{8}b_1(0)+\frac{b_0'(0)}{2}=\frac{b_1(0)}{8}=\frac{1}{16}u_\alpha^+(P_+)\,$$
due to (\ref{rel2}) and (\ref{rel3}). All the $a$-periods of this
differential vanish, therefore, it coincides with
$$\frac{1}{16}u_\alpha^+(P_+)W_+(\,\cdot\,, P_+)$$
and the first asymptotics in (\ref{adif1}) is proved.

The other asymptotics of Theorem \ref{Abdiff} can be proved in a similar way.
Theorem \ref{bergas} follows from Theorem \ref{Abdiff} (see \cite{Fay} p. 41
for a short explanation of this implication). $\square$

It is also possible to prove Theorem \ref{bergas} independently: one starts
from the generalization of Lemma 1 given in Remark \ref{Polusa}, using this
generalization with, say, $v_-=0$ and $v_+=W_+(\,\cdot\,, Q)$ with $Q\in
\L^+\setminus D^+$, one establishes expansion (\ref{MAIN}) for one-form
$W(\,\cdot\,, Q)$ exactly in the same manner as it was done for a basic
holomorphic differential. Repeating the proof of Theorem \ref{Abdiff} with
$W(\,\cdot\,, Q)$ instead of $v^{(s)}_\alpha$ we arrive to the asymtotics
stated in Theorem \ref{bergas}.

The following proposition gives the asymptotics of  other type than given in
Theorem \ref{bergas}: now one of the arguments of the canonical meromorphic
bidifferential lies inside the pinching zone (being one of the two endpoints
of the cut).
\begin{proposition}\label{collideas}
Let a point $P$ lies on the surface $\L^{\pm}$ far from the pinching zone and
let $P_r=\l^{-1}(s)$ and $P_l=\l^{-1}(0)$ be the critical points of the map
$\l:{\cal U}\to \{\l: |\l|<1\}$. Then

\begin{equation}\label{pravkon}
W(P, P_r)=\frac{\sqrt{s}}{2}W_\pm(P_\pm, P)+O({s^{3/2}}),
\end{equation}
\begin{equation}\label{levkon}
W(P, P_l)=-i\frac{\sqrt{s}}{2}W_\pm(P_\pm, P)+O(s^{3/2}),
\end{equation}
as $s\to 0$. Here the differentials are calculated in the local parameters
related to corresponding branched coverings: i. e. $\sqrt{\l(\cdot)-s}$ at
$P_r$, $\sqrt{\l(\cdot)}$ at $P_l$; $\l^\pm(\cdot)$ at  $P_\pm$  and an
arbitrary local parameter inherited from $\L^\pm$ at $P$.
\end{proposition}

{Proof.} For the $1$-form $W(\,\cdot\,, P)$ one has the expansion
(\ref{MAIN}) with $b_0(0)=0$, $b'_0(0)+\frac{1}{2}b_1(0)=0$ and
$b_1(0)=a_0(0)=\frac{1}{2}W_\pm(P_\pm, P)$.

Now substituting in this expansion $\l=s+t^2$, $d\l=2t\,dt$ setting $t=0$ and
then sending $s\to 0$ we get (\ref{pravkon}). Substituting $\l=t^2$,
$d\l=2t\,dt$,  setting $t=0$ and sending $s\to 0$, we get (\ref{levkon}).
$\square$

\subsection{Asymptotics of $E(P, Q), \sigma (P, Q)$, $C(P)$ and $K^P$}

First recall the following expression, relating the prime form, $E(x, y)$, to the canonical
meromorphic differential on an arbitrary compact Riemann surface of genus $g$
(see \cite{Fay}, p. 26):
\begin{equation}\label{P-C}
\frac{\theta(\int_x^y\vec{v}-e)\theta(\int_x^y\vec{v}+e)}{\theta^2(e) E^2(x,
y)}=W(x, y)+\sum_{i, j=1}^g\frac{\partial^2\log \theta(e)}{\partial
z_i\partial z_j}v_i(x)v_j(y),
\end{equation}
where $\vec{v}=(v_1, \dots, v_g)^t$ is a column of basic holomorphic
differentials, $e$ is an arbitrary vector from ${\mathbb C}^n$.

From this expression taken together with the asymptotics for the basic
holomorphic differentials and the canonical meromorphic bidifferential one
easily derives the following asymptotics for the prime form on the family
$\L_s$.
\begin{itemize}
\item

\begin{equation}\label{pr1}
E^2(P, Q)=E_\pm^2(P, Q)+o(1)
\end{equation}
as $s\to 0$, here the points $P, Q$ belong to $\L^\pm$ and are far from the
pinching zone, $E_\pm(P, Q)$ is the prime form on $\L^\pm$, all the prime
forms are calculated in local parameters near $P$ and $Q$ inherited from
$\L^\pm$;

\item
\begin{equation}\label{pr2}
E^2(P, Q)=-\frac{16}{s^2}E_\pm^2(P, P_\pm)E_\mp^2(Q, P_\mp)+O(\frac{1}{s})
\end{equation}
if $P\in \L^\pm$ and $Q\in \L^\mp$;
\item
\begin{equation}\label{pr3}
E^2(P, P_r)=\frac{2}{\sqrt{s}}E^2_\pm(P, P_\pm)+O(\sqrt{s}),\ \ \ \ \ E^2(P,
P_l)=\frac{2i}{\sqrt{s}}E^2_\pm(P, P_\pm)+O(\sqrt{s}),
\end{equation}
if $P\in \L^\pm$, the local parameter at $P_l$ is $\sqrt{\l}$, the local
parameter near $P_r$ is $\sqrt{\l-s}$.
\end{itemize}

From now on we use the following notation $\Delta(s)=\frac{4}{s}$ and denote
by a single letter $\epsilon$  different unitary constants ("phase factors",
($|\epsilon|=1$) which may appear as additional factors in some of our
formulas; the concrete values of these factors are of no interest for us.

The next two quantities whose asymptotics we need are defined as follows (see
\cite{Fay92}, (1.13) and (1.17)):
\begin{equation}\label{sigma1}
\sigma(P, Q)=\exp\left\{-\sum_{k=1}^g\int_{a_k}v_k(x)\log\frac{E(x, P)}{E(x,
Q)}\right \}\,,\end{equation} and
\begin{equation}\label{defC}
C(P)=\frac{\theta(\int_{P}^{Q_1}\vec{v}\dots
+\int_{P}^{Q_g}\vec{v}+K_P)\prod_{i<j}^gE(Q_i, Q_j)\prod_{i=1}^g\sigma(Q_i,
P)}{ {\rm det}(v_i(Q_j))\prod_{i=1}^gE(P, Q_i)}\,
\end{equation}
where $Q_1, \dots Q_g$ are arbitrary points of $\L$ (expression (\ref{defC})
is independent of the choice of these points) and $K_P$ is the vector of
Riemann constants.

 Using asymptotics for the prime-form
(\ref{pr1}--\ref{pr3}) and the basic holomorphic differentials one easily
obtains from (\ref{sigma1}) the following asymptotics as $s\to 0$:
\begin{equation}\label{assigm1}
\sigma(P, Q)\sim\sigma_\pm(P, Q)\left[\frac{E_\pm(Q, P_\pm)}{E_\pm(P, P_\pm)}
\right]^{g^\mp},
\end{equation}
for $P, Q\in \L^\pm$;
\begin{equation}\label{assigm2}
\sigma(P, Q)\sim\epsilon \sigma_\pm(P, P_\pm)\sigma_\mp(P_\mp,
Q)\frac{\left[E_\mp(P_\mp, Q) \right]^{g^\pm}}{\left[E_\pm(P,
P_\pm)\right]^{g^\mp}}\left[\Delta(s)\right]^{g^\pm-g^\mp}\,,
\end{equation}
if $P\in \L^\pm$, $Q\in \L^\mp$;
\begin{equation}\label{assigm3}
\sigma(P_r, Q)\sim \epsilon \sigma(P_l, Q)\sim \epsilon \sigma_\pm(P_\pm,
Q)\left[E_\pm(P_\pm,
Q)\right]^{g^\mp}\left[\Delta(s)\right]^{(3g^\mp-g^\pm)/4}\,
\end{equation}
if $Q\in \L^\pm$.

The asymptotics of (\ref{defC}) is a bit more tricky to obtain and we give
more details. First choose the points $\{Q_i\}$ in such a way that $g^+$ of
them, $R_1, \dots, R_{g^+}$ belong to $\L^+$ and the other $g^-$ points,
$S_1, \dots, S_{g^-}$, belong to $\L^-$. Then, assuming for definiteness
$P\in \L^+$, one has as $s\to 0$
$$\theta(\int_{P}^{Q_1}\vec{v}\dots +\int_{P}^{Q_g}\vec{v}+K_P| {\mathbb B})\sim$$
$$\theta\left(\int_P^{R_1}\left(
^{\vec{v}_+}_{\vec{0}}\right)+\dots+\int_P^{R_{g^+}}\left(
^{\vec{v}_+}_{\vec{0}}\right)+g_-\int_P^{P_+}\left(^{\vec{v}_+}_{0}\right)+\int_{P_-}^{S_1}\left(^{\vec{0}}_{\vec{v}_-}\right)+\dots
+\int_{P_-}^{{S_{g^-}}}\left(^{\vec{0}}_{\vec{v}_-}\right)+\right.$$$$\left.+\left(^{K_P^+-g^-\int_P^{P_+}\vec{v}_+}_{K_{P_-}^-}
\right)\Big|{\rm diag}({\mathbb B}^+, {\mathbb B}^-)\right)=
$$
\begin{equation}\label{theta1}
=\theta_+(\int_P^{R_1}\vec{v}_++\dots+\int_P^{R_{g^+}}\vec{v}_++K_P^+)\,
\,\theta_-(\int_{P_-}^{S_1}\vec{v}_-+\dots+\int_{P_-}^{S_{g^-}}\vec{v}_-+K_{P_-}^-)\,.
\end{equation}
Now using the asymptotics for the prime form and $\sigma$, we see that the
numerator of (\ref{defC}) (with the just made choice of $Q_1, \dots, Q_g$) is
equivalent to
$$\epsilon\,\theta_+(\int_P^{R_1}\vec{v}_++\dots+\int_P^{R_{g^+}}\vec{v}_++K_P^+)\,
\,\theta_-(\int_{P_-}^{S_1}\vec{v}_-+\dots+\int_{P_-}^{S_{g^-}}\vec{v}_-+K_{P_-}^-)\,\prod_{i<j}E_+(R_i,
R_j)\prod_{i<j}E_-(S_i, S_j)
$$
$$\left\{\prod_{i=1}^{g^+}\prod_{j=1}^{g^-}E_+(R_i,
P_+)E_-(S_j,P_-)\right\}[\Delta(s)]^{g^+g^-} \prod_{i=1}^{g^+}\sigma_+(R_i,
P) \frac{ \left\{ E_+(P, P_+)\right\}^{g^+g^-}} {
\left\{\prod_{j=1}^{g^+}E_+(R_j, P_+) \right\}^{g^-} }
$$
$$[\sigma_+(P_+, P)]^{g^-}[E_+(P_+,
P)]^{(g^-)^2}[\Delta(s)]^{g^-(g^--g^+)}\prod_{j=1}^{g^-}\frac{\sigma_-(S_j,
P_-)}{\{E_-(S_j, P_-)\}^{g^+}}\,,
$$
whereas the denominator of (\ref{defC}) is equivalent to
$$\epsilon\left\{ \prod_{i=1}^{g^+}E_+(P, R_i)\right\} [E_+(P, P_+)]^{g^-}\left\{\prod_{j=1}^{g^-}E_-(P_-, S_j)\right\}[\Delta(s)]^{g^-}{\rm det}(v^+_i(R_j))
{\rm det}(v^-_i(S_j))\,.$$ So, after rearranging the terms and numerous
cancelations, one gets the asymptotics
\begin{equation}\label{asC}
C(P)\sim \epsilon C_\pm(P)C_\mp(P_\mp)\left\{E_\pm(P,
P_\pm)\right\}^{g^\mp(g^\pm+g^\mp-1)}\left\{\sigma_\pm(P_\pm,
P)\right\}^{g_\mp}\Delta(s)^{[g^\mp]^2-g^\mp}\,
\end{equation}
if $P\in\L^\pm$.

\begin{remark} {\rm Let us emphasize that in order to define the vector $K^P$ and the
Abel map ${\cal A}_P$ (as well as the prime-form and the left hand side of
expression (\ref{P-C})) one has to introduce the system of cuts on the
surface $\L$ in such a way that the integration $\int_x^y \vec{v}$ is
well-defined for any $x, y$ belonging to the surface $\L=\L_s$ dissected
along the cuts. We choose  the usual symplectic  basis   of homologies
$\{a_\alpha^\pm, b_\alpha^\pm\}_{\a=1, \dots, g^\pm}$ on $\L^\pm$, take
curves representing this basis and dissect the $\L^\pm$ along these curves.
The resulting dissected surface $\L^\pm$ is homeomorphic to a sphere with
$g^\pm$ holes, whereas the surface $\L_s$ dissected along the same curves is
homeomorphic to a sphere with $g$ holes. Notice that the boundary of any hole
is the trivial cycle ($a_\alpha^\pm+ b_\alpha^\pm
-a_\alpha^\pm-b_\alpha^\pm=0$) and, therefore, the $\int_x^y \vec{v_\pm}$ and
$\int_x^y \vec{v}$ are well-defined on the corresponding dissected surfaces.}
\end{remark}
The following lemma immediately follows from the definition of
the vector of the Riemann constants,
$$K^{P}_\beta=\frac{1}{2}+\frac{{\mathbb B}_{\beta\beta}}{2}-\sum_{\alpha=1, \alpha\neq
\beta}^g\int_{a_\alpha}\left(v_\alpha\int_P^xv_\beta\right),$$ and Theorem
\ref{Abdiff}.
\begin{lemma}\label{askp}
One has the asymptotics
\begin{equation}\label{asRiem}
K^P\sim \left(^{K_+^P-g^-\int_P^{P_+}\vec{v}_+} _{K^{P_-}_-} \right)\ \ ,
\end{equation}
as $s\to 0$, where $K^P_+$ and $K^{P_-}_-$ are the vectors of Riemann
constants for the surfaces $\L^+$ and $\L^-$ with the base points $P$ and
$P_-$ respectively.
\end{lemma}

\subsection{Asymptotics of $\tau_g$}

Now we are able to prove asymptotics (\ref{astau}) from the Introduction.
Let $M_\pm=2g^\pm-2$ and let $(\omega^\pm)=\sum_{k=1}^{M_\pm}D_k^\pm$ be the divisor of the holomorphic differential $\omega^\pm$ on $\X^\pm$.

Assume that the point $P$ lies on the component $\L^+$. Using Lemma \ref{askp}, one
can pass to the limit $s\to 0$ in the equation (\ref{rdef}). This results in
the relations
\begin{equation}\label{plus}
{\cal A}^+_{P}((\omega^+))=2K^P_++{\mathbb B}^+{\bf r}^++{\bf q}^+
\end{equation}
and
\begin{equation}\label{minus}
{\cal A}^-_{P_-}((\omega^-))=2K^{P_-}_-+{\mathbb B}^-{\bf r}^-+{\bf q}^+\,,
\end{equation}
where ${\bf r}=({\bf r}^+, {\bf r}^-)$, ${\bf q}=({\bf q}^+, {\bf q}^-)$ and
${\cal A}^\pm$ is the Abel map on $\L^\pm$.

Now one has
$$\tau^{-6}(\L, \l)\sim\epsilon e^{2\pi i <{\bf r}^+, \ K^P_+>}e^{2\pi i <{\bf r}^-,\  K^{P_-}_->}e^{-2\pi
ig^-<{\bf r}^+, \ \int^{P_+}_P\vec{v}_+>}$$$$
\{C_+(P)\}^{-4}\{C_-(P_-)\}^{-4}\{E_+(P, P_+)\}^{4g^-(1-g)}\{\sigma_+(P_+,
P)\}^{-4g^-}[\Delta(s)]^{4(g^--(g^-)^2)}$$
$$\prod_{k=1}^{M_+}\sigma_+(D_k^+, P)\left[\frac{E_+(P, P_+)}{E_+(D_k^+,
P_+)}\right]^{g^-}$$
$$ \left[\sigma_+(P_+, P)\{E_+(P_+,
P)\}^{g^-}[\Delta(s)]^{(3g^--g^+)/4}\right]^2
$$
$$\prod_{k=1}^{M_-}\left\{\sigma_-(D_k^-, P_-)\sigma_+(P_+, P)\frac{[E_+(P_+, P)]^{g^-}}{[E_-(D_k^-,
P_-)]^{g^+}}[\Delta(s)]^{g^--g^+}\right\}$$
$$\prod_{k=1}^{M_+}\{E_+(D_k^+, P)\}^{(g-1)}\left[[\Delta(s)]^{1/4}E(P,
P_+)\right]^{2(g-1)}\prod_{k=1}^{M_-}\{\Delta(s)E_+(P, P_+)E_-(D_k^-,
P_-)\}^{(g-1)}\,$$ with
$g=g^++g^-$. Observe that $\Delta(s)$ enters the above expression with power
$$4(g^--(g^-)^2)+\frac{3g^--g^+}{2}+(g^--g^+)(2g^--2)+\frac{g-1}{2}+(g-1)(2g^--2)=\frac{3}{2}\,,$$
all the factors $E_+(P, P_+)$ cancel out
($4g^-(1-g)+g^-(2g^+-2)+2g^-+g^-(2g^--2)+2(g-1)+(g-1)(2g^--2)=0$) and the
remaining terms can be rearranged into the product of
\begin{equation}\label{tauplus} e^{2\pi i <{\bf r}^+, \
K^P_+>}C_+^{-4}(P)\prod_{k=1}^{M_+}\sigma_+(D_k^+,
P)\left\{E_+(D_k^+, P) \right\}^{(g^+-1)}\,,\end{equation}
\begin{equation}\label{tauminus}
e^{2\pi i <{\bf r}^-,\
K^{P_-}_->}C_-^{-4}(P_-)\prod_{k=1}^{M_-}\sigma_-(D_k^-,
P_-)\left\{E_-(D_k^-, P_-) \right\}^{(g^--1)}\,.\end{equation}
and \begin{equation}\label{tri} e^{-2\pi ig^-<{\bf r}^+, \
\int^{P_+}_P\vec{v}_+>}\left\{[\sigma_+(P_+,
P)]^{-2}\frac{\prod_{k=1}^{M_+}E_+(D_k^+,
P)}{\prod_{k=1}^{M_+}E_+(D_k^+, P_+)} \right\}^{g^-}
\,.\end{equation}
 According to \cite{DG} (see Theorem 2 on page 47), the
 expression in the braces in (\ref{tri}) is nothing but
 $e^{2\pi i<{\bf r}^+,\  {\cal A}^+_{P}(P_+)>}$ and, therefore, the
 expression (\ref{tri}) equals one; expressions (\ref{tauplus}) and
 (\ref{tauminus}) coincide with $\tau_{g^+}^{-6}(\L^+, \omega^+, \{a_\a^+, b_\a^+\})$ and $\tau_{g^-}^{-6}(\L^-, \omega^-, \{a_\a^-, b_\a^-\} )$
 respectively. $\square$

\section{Surgery and asymptotics}
\subsection{Wentworth lemma}
The following important Lemma essentially coincides with the statement proved in \S 3 of \cite{Wentworth}. We formulate the Wentworth result, adapting it for our needs.
\begin{lemma}\label{Wen} Let $\X$ be a translation surface, $\Delta$ the Friedrichs extension of the Laplacian on $\X$. Let $z$ be a local parameter near a (nonsingular) point $P\in \X$ such that $z(P)=0$ and $\Delta=-4\partial_z\partial_{\bar z}$ in the unit ball $B(1)=\{|z|\leq 1\}\subset \X$.
For $0<\epsilon\leq 1$ set $B(\epsilon)=\{|z|\leq \epsilon\}$ and $\X_\epsilon=\X\setminus B(\epsilon)$.
Denote  by ${\cal N}_\epsilon$ the Dirichlet-to-Neumann operator for $\X_\epsilon$:
$${\cal N}_\epsilon: C^\infty(\partial \X_\epsilon)\to C^\infty(\partial \X_\epsilon)\,,$$
$${\cal N}_\epsilon (f)=\partial_nu\Big|_{\partial \X_\epsilon},$$
where  the function $u$ satisfies
\begin{equation}\label{harmo}
\begin{cases}\Delta u=0\ \ \ \ {\text in}\ \ \X_\epsilon\\
u\Big|_{\partial X_\epsilon}=f
\end{cases}
\end{equation}
and $n$ is the unit outer normal to $\partial \X_\epsilon$. (Actually, ${\cal N}_\epsilon$ is a pseudodifferential operator of order $1$ on $\partial \X_{\epsilon}$).

Let $z=re^{i\phi}$, then $\phi$ is the angular coordinate on the circle $\{r=\epsilon\}=\partial X_\epsilon$. Let $f\in L_2(\partial X_\epsilon, d\phi)$, $f(\phi)=\sum_{k\in {\mathbb Z}}a_ke^{ik\phi}$. Define the (unbounded) operators ${\bf \nu}$ and $|{\bf \nu}|$ in $L_2(\partial X_\epsilon, d\phi)$ via
$${\bf \nu}f(\phi)=\sum_{k\in {\mathbb Z}}ka_ke^{ik\phi}$$ and
$$|{\bf \nu}|f(\phi)=\sum_{k\in {\mathbb Z}}|k|a_ke^{ik\phi}\,.$$
Then one has the following relation:
\begin{equation}
\epsilon{\cal N}_\epsilon=|{\bf \nu}|+{\bf O}(\epsilon)\,
\end{equation}
where ${\bf O}(\epsilon)$ is the operator of trace class in $L_2(\partial \X_\epsilon, d\phi)$ with the trace norm which is asymptotically $O(\epsilon)$ as $\epsilon \to 0$.

\end{lemma}

For completeness we give the proof here (it differs from the one given in \cite{Wentworth} by insignificant changes).
Introduce the operator $R_\epsilon:L_2(S^1, d\phi)\to L_2(S^1, d\phi)$ via $R_\epsilon f=g=u|_{\{|z|=1\}}$, where $u$ and $f$ are from (\ref{harmo}). Using Green formula for the Friedrichs extension of the Laplacian, it is easy to check the identity
 $$\int_{|z|=\epsilon}|f|^2d\phi-\int_{|z|=1}|g|^2d\phi=2\int_{\epsilon}^1\frac{dr}{r}\iint_{\X_r}|\nabla u|^2,$$
which implies the norm estimate
\begin{equation}\label{norm}||R_\epsilon||\leq 1\,. \end{equation}
(It is important here that $\Delta$ is the Friedrichs extension; for other extensions the above double integral may be infinite!)

The function $u$ from (\ref{harmo}) is harmonic in the annulus $\{\epsilon\leq|z|\leq 1\}$ and, therefore, admits there the standard representation
$$u(r, \phi)=\sum_{k\in {\mathbb Z}}a_ke^{ik\phi}r^k+c_0\log r+\sum_{k\in {\mathbb Z}\setminus 0}(r^{|k|}-r^{-|k|})b_ke^{ik\phi}\,.$$
On the other hand the Green formula implies the relation $\int_{\partial \X_\epsilon}\frac{\partial u}{\partial r}=0$ and, therefore, one has $c_0=0$ in the previous representation.
Thus, the function $u$ from (\ref{harmo}) is representable inside the annulus $\{\epsilon\leq|z|\leq 1\}$ as
\begin{equation}\label{predstav}
u=a_0+\sum_{k\neq 0}a_kr^ke^{ik\phi}+\sum_{k\neq 0}b_kr^{-k}e^{ik\phi}\,.
\end{equation}
Now notice that the operators $\epsilon {\cal N}_\epsilon$ and $R_\epsilon$ map the boundary value of the function $u$ from (\ref{harmo}) and (\ref{predstav}) at the circle $|z|=\epsilon$ to the functions
$$\sum_{k\neq 0}k(b_k\epsilon^{-k}-a_k\epsilon^k)e^{ik\phi}$$
and
$$a_0+\sum_{k\neq 0}(a_k+b_k)e^{ik\phi}$$
respectively.
For a sequence of complexe numbers $\{\alpha_k\}_{k\in {\mathbb Z}}$
introduce the operator (may be unbounded) ${\rm Op}(\alpha_k)$ in $L_2(S^1, d\phi)$ via
$${\rm Op}(\alpha_k)f (\phi)=\sum_{k\in {\mathbb Z}}\alpha_ka_ke^{ik\phi}\,.$$
where
 $f\in L_2(S^1, d\phi)$, $u=\sum_{k\in {\mathbb Z}}a_ke^{ik\phi}$. (In this notation ${\bf \nu}={\rm Op}(k)$ and $|{\bf \nu}|={\rm Op}(|k|)$.)

Now one has
$${\rm Op}(\epsilon^k-\epsilon^{-k})(\epsilon {\cal N}_\epsilon) (u|_{|z|=\epsilon})=\sum_{k\neq 0}k(b_k+a_k-a_k\epsilon^{2k}-b_k\epsilon^{-2k})e^{ik\phi}=$$
$$=\sum_{k\neq 0}(ka_k+kb_k)e^{ik\phi}-\sum_{k\neq o}\{k(a_k\epsilon^k+b_k\epsilon{-k})(\epsilon^k+\epsilon^{-k})-k(a_k+b_k)\}e^{ik\phi}=$$
$$=2\sum_{k\neq 0}(ka_k+kb_k)e^{ik\phi}-\sum_{k\neq 0}(\epsilon^k+\epsilon^{-k})k(a_k\epsilon^k+b_k\epsilon^{-k})e^{ik\phi}\,$$
which implies the relation
\begin{equation}\label{key}
\epsilon{\cal N}_\epsilon={\rm Op}(\frac{2}{\epsilon^k-\epsilon^{-k}}){\bf \nu}R_\epsilon-{\rm Op}(\frac{\epsilon^k+\epsilon^{-k}}{\epsilon^k-\epsilon^{-k}}){\bf \nu}\,.
\end{equation}

(Notice that  functions from the image of the operator ${\bf \nu}$ are orthogonal to $1$ and, therefore,  the  right hand side of (\ref{key}) is correctly defined.)
Clearly, the operator ${\rm Op}(\frac{2}{\epsilon^k-\epsilon^{-k}}){\bf \nu}$ is of trace class with the trace norm $|||{\rm Op}(\frac{2}{\epsilon^k-\epsilon^{-k}}){\bf \nu}|||=O(\epsilon)$, due to (\ref{norm}) the same is true for the first term in the right hand side of (\ref{key}). For $k\neq 0$ one has $\frac{\epsilon^k+\epsilon^{-k}}{\epsilon^k-\epsilon^{-k}}\to -{\text sgn}\,k$ as $\epsilon \to 0$
and the simple estimate shows that
$$-{\rm Op}(\frac{\epsilon^k+\epsilon^{-k}}{\epsilon^k-\epsilon^{-k}}){\bf \nu}=|{\bf \nu}|+ r(\epsilon)\,,$$
where $|||r(\epsilon)|||=O(\epsilon^2)$ which proves the Lemma.
\subsection{Analytic surgery for translation surfaces}
The following proposition is a variant of Theorem B* from \cite{BFK}. Its proof does not essentially differ from the proof of classical BFK formula.
\begin{proposition}\label{BFK}
Let $\Gamma$ be a smooth closed curve on a translation surface $\X$ containing no conical points and dividing $\X$ into two parts $\X_1$ and $\X_2$ with common boundary $\Gamma$. Let $(\Delta, \X_{1,2})$ be the operators of the Dirichlet boundary value problems in $\X_{1, 2}$.
Then one has the relation
$${\rm det}^* \Delta=\frac{{\rm Area}(X)}{{\rm length}(\Gamma)}{\rm det}(\Delta, \X_1){\rm det}(\Delta, \X_2){\rm det}^*({\cal N}_1+{\cal N}_2),                                $$
where ${\cal N}_{k}$ is the Dirichlet-to-Neumann operator $C^{\infty}(\Gamma)\to \C^\infty(\Gamma)$,
${\cal N}_{k}(f)=\partial_{n_{k}}u_{k}|_\Gamma$ with $\Delta u_k=0$ in $\X_{k}$, $u|_\Gamma=f$ and $n_{k}$ being the outer unit normal to $\partial \X_{k}$, $k=1,2$.
\end{proposition}

The following proposition (see \cite{Wentworth}) is a consequence of Wentworth lemma
\begin{proposition}
Let $\X_1=\X\setminus B(\epsilon)$, $\X_2=B(\epsilon)$, $\Gamma=\{|z|=\epsilon\}$. Then
\begin{equation}\label{lim}\lim_{\epsilon \to 0}\frac{ {\rm det}^*({\cal N}_1+{\cal N}_2)}{{\rm length}\, (\Gamma)}=\frac{1}{2}\,.\end{equation}
\end{proposition}
We give a proof of this proposition following \cite{Wentworth}.
Representing the function $u_2$ harmonic in the disk $\X_2=\{|z|\leq \epsilon\}$  in the form
$$u_2=\sum_{k\in {\mathbb Z}}a_kr^{|k|}e^{ik\phi}\,,$$
one immediately gets the relation
$$\epsilon{\cal N}_2=|{\bf \nu}|\,.$$
One has now
\begin{equation}\label{one}
\epsilon({\cal N}_1+{\cal N}_2)=2|{\bf \nu}|+{\bf O}(\epsilon)
\end{equation}
and, therefore,
\begin{equation}\label{two}
\log {\rm det}^*\{\epsilon({\cal N}_1+{\cal N}_2)\}=\log {\rm det}^*\{2|{\bf \nu}|\}+o(1)
\end{equation}
which implies
\begin{equation}\label{three}
\log\frac{{\rm det}^*({\cal N}_1+{\cal N}_2)}{\epsilon}=\log \frac{{\rm det}^*|{\bf \nu}|}{2}+o(1)
\end{equation}
or, what is the same,
\begin{equation}\label{four}
\log\frac{{\rm det}^*({\cal N}_1+{\cal N}_2)}{2\pi\epsilon}=-\log 2+\log {\rm det}^*|\nu|-\log{2\pi}+o(1).
\end{equation}
Using the known properties of the Riemann zeta-function,  $-2\zeta'(0)=\log{2\pi}$ and $\zeta(0)=-1/2$, one gets the relation
${\rm det}^*|\nu|=2\pi$ which (together with (\ref{four})) implies the (\ref{lim}).
\begin{remark}{\rm
Implication (\ref{one})$\Rightarrow$(\ref{two}) is a consequence of the following estimate
$$|\log {\rm det}^*(2|\nu|+{\bf O}(\epsilon))-\log {\rm det}^*(2|\nu|)|=\left|\int_0^1\frac{d}{dt}\log {\rm det}^*(2|\nu|+t{\bf O}(\epsilon))dt   \right|\leq$$
$$\leq\int_0^1\left|{\rm tr} \left((2|\nu|+t{\bf O}(\epsilon))\Big|_{\{1\}^\bot}^{-1}{\bf O}(\epsilon)\right)\right|dt\leq C_1|{\rm tr}{\bf O}(\epsilon)|\leq C_2\epsilon $$
(cf. \cite{Lee}, Lemma 4.1).}

\end{remark}
\begin{remark}\label{zetaotnulja}    {\rm
Implication (\ref{two})$\Rightarrow$(\ref{three}) is a consequence
of the standard relations
\begin{equation}\label{rescale}{\rm det}^*(\epsilon A)=\epsilon^{\zeta_A(0)}{\rm det}^*A,\end{equation}
\begin{equation}\label{heat1}
\zeta_{{\cal N}_1+{\cal N}_2}(0)={\bf h}_0-{\rm dim}\,{\rm Ker}({\cal N}_1+{\cal N}_2)={\bf h}_0-1,
\end{equation}
where ${\bf h}_0$ is the constant term in the asymptotical expansion of ${\rm tr}\,e^{-t({\cal N}_1+{\cal N}_2)}$ as $t\to 0+$,
the relation $\zeta_{|\nu|}(0)=-1$
and the result from \cite{EdwardWu}:
\begin{equation}\label{nul}{\bf h}_0=0\,.\end{equation}

}
\end{remark}

\begin{corollary}
One has the asymptotics
\begin{equation}\label{vyrezdisk}
{\rm det}(\Delta, \X\setminus B(\epsilon))\sim\frac{2^{7/6}\sqrt{\pi}e^{2\zeta'(-1)+5/12}{\rm det}^*\Delta}{{\rm Area}(\X)}\epsilon^{1/3}
\end{equation}
as $\epsilon\to 0$.
\end{corollary}
{\bf Proof}. This immediately follows from Proposition \ref{BFK}, (\ref{lim}) and the relation
$${\rm det}(\Delta, B(\epsilon))=2^{-1/6}\pi^{-1/2}\epsilon^{-1/3}e^{-2\zeta'(-1)-5/12}$$
found in \cite{Weisberger}.

\begin{remark}{\rm It is interesting to compare (\ref{vyrezdisk}) with the asymptotics of the first  eigenvalue of the operator of the Dirichlet boundary value problem in $\X\setminus B(\epsilon)$,
$$\lambda_1(\Delta, \X\setminus B(\epsilon))\sim -\frac{2\pi}{{\rm Area}(\X)}(\log \epsilon)^{-1}$$
as $\epsilon \to 0$, which was found in
\cite{Ozawa}.}
\end{remark}

\subsection{Symmetric case}
Let $\X$ be a translation surface of genus $g\geq 1$, $z$ a local coordinate in a vicinity of a nonsingular point $P$ of $\X$, such that $z(P)=0$ and in the ball
$\{|z|\leq \epsilon\}$ the operator $\Delta$ acts as $-4\partial_z\partial_{\bar z}$. Introduce a straight cut $I(\epsilon/2)$ connecting the points $z=0$ and $z=\epsilon/2$ and glue two copies of $\X\setminus I(\epsilon/2)$ along the cut in a usual way. One gets a translation surface $\widehat{\X}$ of genus  $2g$ and the area ${\rm Area}(\widehat\X)=2{\rm Area}(\X)$. The end points of the cut give rise to two conical points, $P_1, P_2$ of conical angles $4\pi$ on $\widehat{\X}$.
Let $\widehat \Delta$ be the (Friedrichs extension of)  Laplacian on $\widehat \X$.
The following statement is a very special case of (\ref{Main1}) proved in an alternative way in order to get information about the unknown constants $\delta_g$ in (\ref{Main1})
\begin{proposition}
One has the asymptotics
\begin{equation}\label{AuxAs}
{\rm det}^*\widehat \Delta\sim\frac{2\kappa_0}{{\rm Area}(\X)}\left\{{\rm det}^*\Delta\right\}^2\epsilon^{1/2},
\end{equation}
as $\epsilon \to 0$, where the  constant $\kappa_0$ is the same for all translation surfaces $\X$  (and for all $g\geq 1$) and is defined via formula (\ref{kappa0}) below.
\end{proposition}
\begin{remark} {\rm
The factor $2/{\rm Area} (\X)$ is nothing but ${\rm Area} (\widehat \X)/[{\rm Area}(\X){\rm Area} (\X)]$ that is why we are not attaching the factor $2$ to the constant $\kappa_0$ in (\ref{AuxAs}).}
\end{remark}
{\bf Proof.} First, notice that the surface $\widehat \X$ is provided with a natural involution $*$, and the shores of the cut $I(\epsilon)$ (two homologous saddle connections, $\gamma$ and $\gamma'$ on $\widehat \X$) are fixed by this involution.
One has the standard ($\widehat \Delta$-invariant) decomposition $L_2(\widehat \X)=L_2^{symm}(\widehat \X)\bigoplus L_2^{antisymm}(\widehat \X)$
and the functions $u$ from the domain of $\widehat \Delta$ which enter $L_2^{antisymm}(\widehat \X)$ satisfy $u|_{\gamma\cup \gamma'}=0$, whereas
the functions $u$ from the domain of $\widehat \Delta$ which enter $L_2^{symm}(\widehat \X)$ satisfy $u_n|_{\gamma\cup \gamma'}=0$.
This shows that the operator $\widehat \Delta$ is unitary equivalent to the direct sum of the operators $\Delta_{\cal D}$ and $\Delta_{\cal N}$ of the homogeneous Dirichlet and Neumann boundary value problems in $\X\setminus I(\epsilon)$ (cf., e. g., \cite{Hillairet}, p. 79) and, therefore,
\begin{equation}\label{proizv}{\rm det}^*\widehat\Delta={\rm det}\Delta_{\cal D}{\rm det}^*\Delta_{\cal N}\,.\end{equation}
We are to study the asymptotics of ${\rm det}\Delta_{\cal D}$ and ${\rm det}^*\Delta_{\cal N}$ as $\epsilon\to 0$.

{\bf Asymptotics of $ {\rm det}\Delta_{\cal D}$.}  We will be using the generalizations of the BFK formula (Theorems B and B* from \cite{BFK}) to the case of Laplacians on $2d$ manifolds with boundary with Dirichlet (and Neumann) boundary conditions. Such generalizations are straightforward and are mentioned in \cite{Lee} (see Remark on page 326). Their proofs differ from the standard proof of Theorem B* from \cite{BFK} insignificantly.
For the operator $\Delta_{\cal D}$ the following surgery formula holds true:
\begin{equation}\label{SurDir}
{\rm det}\Delta_{\cal D}={\rm det}(\Delta, \X\setminus B(\epsilon)){\rm det}(\Delta, B(\epsilon)\setminus I(\epsilon)){\rm det}({\cal N}_\epsilon+{\cal N}^{int, D}_\epsilon)\,,
\end{equation}
where ${\cal N}_\epsilon$ is the Dirichlet-to-Neumann operator from Lemma \ref{Wen}; the operator ${\cal N}_\epsilon^{int, D}:C^\infty (\partial B(\epsilon))\to C^\infty(\partial B(\epsilon))$ is defined via ${\cal N}^{int, D}_\epsilon(f)=u_n|_{\partial B(\epsilon)}$, with $u$ subject to
\begin{equation}\label{dir1}\begin{cases}
\Delta u=0 \ \ {\rm in}\ \ B(\epsilon)\setminus I(\epsilon)\\
u|_{\partial B(\epsilon)}=f\\
u|_{I(\epsilon)}=0
\end{cases}\end{equation}
and $(\Delta, B(\epsilon)\setminus I(\epsilon))$ is the operator of the homogeneous Dirichlet boundary value problem in $B(\epsilon)\setminus I(\epsilon)$.
(Notice that there are no coefficient of the type ${\rm Area}/{\rm length}$ at the right hand side of (\ref{SurDir}): all the operators there are invertible and (\ref{SurDir}) is an analog of Theorem B from \cite{BFK}.) The asymptotics of the first factor in (\ref{SurDir}) is given in (\ref{vyrezdisk}), the asymptotics of other two factors can be obtained as consequences of homogeneity properties.
Due to (\ref{rescale}) one has
\begin{equation}\label{resc1}
{\rm det}(\Delta, B(\epsilon)\setminus I(\epsilon))=\epsilon^{-2\zeta_{(\Delta, B(\epsilon)\setminus I(\epsilon))}(0)}{\rm det}(\Delta, B(1)\setminus I(1)),
\end{equation}
where the value of $\zeta_{(\Delta, B(\epsilon)\setminus I(\epsilon))}(0)$ coincides with the term ${\bf h}_0$ of the corresponding heat asymptotics (cf., (\ref{heat1}); clearly, ${\rm dim}{\rm Ker}(\Delta, B(\epsilon)\setminus I(\epsilon))=0$). The term ${\bf h}_0$ is easy to find, namely
one has
\begin{equation}
{\bf h}_0=\frac{1}{6}+2\frac{\pi^2-(2\pi)^2}{24\pi(2\pi)}=\frac{1}{24},
\end{equation}
where the term $\frac{1}{6}$ comes from the  part $\partial B(\epsilon)$ of $\partial[B(\epsilon)\setminus I(\epsilon)]$ and two terms
$\frac{\pi^2-(2\pi)^2}{24\pi(2\pi)}$ come from two angle points of opening $\beta=2\pi$ at the end points of the cut $I(\epsilon)$ (see \cite{Cheeger}, formula (4.41) or \cite{AurellSalomonson}, formula (37); the straight part, ${\rm int}(I(\epsilon))$ of the boundary makes no input in ${\bf h}_0$).
Thus, one has
\begin{equation}\label{vtorterm}
{\rm det}(\Delta, B(\epsilon)\setminus I(\epsilon))=\epsilon^{-\frac{1}{12}}{\rm det}(\Delta, B(1)\setminus I(1))\,.
\end{equation}
Moreover, from the result of \cite{EdwardWu} and the relation ${\rm dim}{\rm Ker}({\cal N}_\epsilon+{\cal N}^{int, D}_\epsilon)=0$ one gets
the equality
$$\zeta_{{\cal N}_\epsilon+{\cal N}^{int, D}_\epsilon}(0)=0$$
and, therefore,
$${\rm det}({\cal N}_\epsilon+{\cal N}^{int, D}_\epsilon)={\rm det}(\epsilon{\cal N}_\epsilon+\epsilon{\cal N}^{int, D}_\epsilon)=
{\rm det}(|\nu|+{\bf O}(\epsilon)+{\cal N}^{int,D}_1)\sim {\rm det}(|\nu|+{\cal N}^{int, D}_1),$$
as $\epsilon \to 0$  due to Lemma \ref{Wen}. Summarizing, one arrives at the asymptotics
\begin{equation}\label{maindir}
{\rm det}\Delta_{\cal D}\sim\frac{2^{7/6}\sqrt{\pi}e^{2\zeta'(-1)+5/12}{\rm det}^*\Delta\,{\rm det}(|\nu|+{\cal N}^{int, D}_1)   {\rm det}(\Delta, B(1)\setminus I(1))}{{\rm Area}(\cal X)}\epsilon^{\frac{1}{4}}
\end{equation}
as $\epsilon\to 0$.

{\bf Asymptotics of $ {\rm det}\Delta_{\cal N}$.}
For the operator $\Delta_{\cal N}$ the analog of the Theorem B* from \cite{BFK} looks as follows:
\begin{equation}\label{SurDir1}
{\rm det}^*\Delta_{\cal N}=\frac{{\rm Area}(\X)}{2\pi \epsilon}{\rm det}(\Delta, \X\setminus B(\epsilon)){\rm det}(\Delta, B(\epsilon)\setminus I(\epsilon); D, N){\rm det}^*({\cal N}_\epsilon+{\cal N}^{int, N}_\epsilon)\,,
\end{equation}
where ${\cal N}_\epsilon$ is the Dirichlet-to-Neumann operator from Lemma \ref{Wen}; the operator ${\cal N}_\epsilon^{int, N}:C^\infty (\partial B(\epsilon))\to C^\infty(\partial B(\epsilon))$ is defined via ${\cal N}^{int, N}_\epsilon(f)=u_n|_{\partial B(\epsilon)}$, with $u$ subject to
\begin{equation}\label{neum1}\begin{cases}
\Delta u=0 \ \ {\rm in}\ \ B(\epsilon)\setminus I(\epsilon)\\
u|_{\partial B(\epsilon)}=f\\
u_n|_{I(\epsilon)}=0
\end{cases}\end{equation}
and $(\Delta, B(\epsilon)\setminus I(\epsilon); D, N)$ is the operator of the homogeneous boundary value problem in $B(\epsilon)\setminus I(\epsilon)$ with Dirichlet conditions on $\partial B(\epsilon)$ and Neumann conditions on $I(\epsilon)$.
As above,  the asymptotics of the first factor in (\ref{neum1}) is given in (\ref{vyrezdisk}), the asymptotics of other two factors can be obtained as consequences of homogeneity properties.
Due to (\ref{rescale}) one has
\begin{equation}\label{resc2}
{\rm det}(\Delta, B(\epsilon)\setminus I(\epsilon))=\epsilon^{-2\zeta_{(\Delta, B(\epsilon)\setminus I(\epsilon); D, N)}(0)}{\rm det}(\Delta, B(1)\setminus I(1); D, N)=\epsilon^{-\frac{1}{12}}{\rm det}(\Delta, B(1)\setminus I(1); D, N)\,.
\end{equation}
(The inputs from the angle points to the ${\bf h}_0$ are the same for Dirichlet and Neumann problems.)
Since  ${\rm dim}{\rm Ker}({\cal N}_\epsilon+{\cal N}^{int, N}_\epsilon)=1$ one gets
the equality
$$\zeta_{{\cal N}_\epsilon+{\cal N}^{int, N}_\epsilon}(0)=-1$$
and, therefore,
$${\rm det}^*({\cal N}_\epsilon+{\cal N}^{int, N}_\epsilon)=\epsilon{\rm det}^*(\epsilon{\cal N}_\epsilon+\epsilon{\cal N}^{int, N}_\epsilon)=
\epsilon{\rm det}^*(|\nu|+{\bf O}(\epsilon)+{\cal N}^{int,N}_1)\sim \epsilon {\rm det}^*(|\nu|+{\cal N}^{int, N}_1),$$
as $\epsilon \to 0$  due to Lemma \ref{Wen}. Summarizing, one arrives at the asymptotics
\begin{equation}\label{mainneum}
{\rm det}^*\Delta_{\cal N}\sim 2^{1/6}\pi^{-1/2}e^{2\zeta'(-1)+5/12}{\rm det}^*\Delta\,{\rm det}^*(|\nu|+{\cal N}^{int, N}_1)   {\rm det}(\Delta, B(1)\setminus I(1); D, N)\epsilon^{\frac{1}{4}}
\end{equation}
as $\epsilon\to 0$.
\begin{figure}[hbt]
\centering
\includegraphics[scale=0.2]{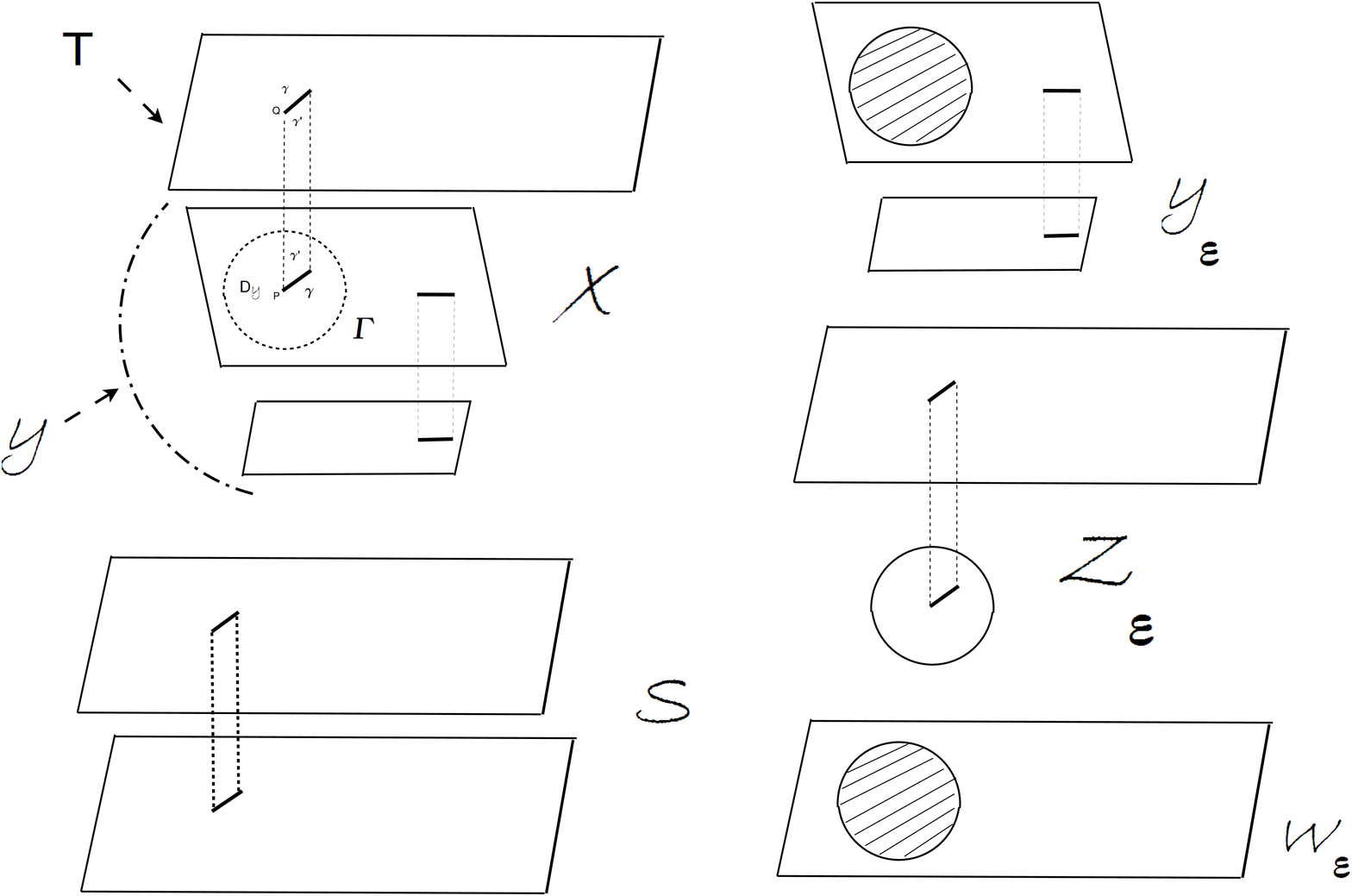}
\caption{}
\la{fig1}
\end{figure}
Now from (\ref{mainneum}), (\ref{maindir}) and (\ref{proizv}) one gets (\ref{AuxAs}) with
\begin{equation}\label{kappa0}
\kappa_0=2^{1/3}e^{4\zeta'(-1)+5/6}{\rm det}(|\nu|+{\cal N}^{int, D}_1)   {\rm det}(\Delta, B(1)\setminus I(1))\,{\rm det}^*(|\nu|+{\cal N}^{int, N}_1)   {\rm det}(\Delta, B(1)\setminus I(1); D, N)\,.
\end{equation}
$\square$

Now from (\ref{Main1}) (with $s=\epsilon/2$) and (\ref{AuxAs}) one gets the relation
$$\delta_{2g}=2\sqrt{2}\kappa_0(\delta_g)^2$$
for the constant $\delta_g$ from (\ref{main}). This implies
\begin{equation}\label{valuedelta}\delta_{N}=(2\sqrt{2}\kappa_0)^{N-1}\delta_1^{N},\end{equation}
with $\delta_1$ from (\ref{delta1}) for any $N$ of the form $N=2^n$. In the next subsection we show that (\ref{valuedelta}) holds for any natural number $N$.
\subsection{General case} Let ${\cal Y}$ be a translation surface of genus $g-1$ and let also $T$ be a translation surface of genus one (a flat torus). Take two disks, $D_{\cal Y}$ and $D_{T}$  of radius $\epsilon$ in ${\cal Y}$ and $T$ with centers $P$ and $Q$ and introduce two straight cuts of length $\epsilon/2$ starting at points $P\in {\cal Y}$ and $Q\in T$. Gluing the surfaces ${\cal Y}$ and $T$ along the cuts one gets the translation surface ${\cal X}$ of genus $g$. The shores of the cuts give rise to the saddle connections $\gamma$ and $\gamma'$ on $\X$. The boundary $\partial B(\epsilon)=\partial D_{{\cal Y}}$ of the disk in ${\cal Y}$ gives rise to the contour $\Gamma$ on $\X$. Let ${\cal Y}_\epsilon={\cal Y}\setminus D_{{\cal Y}}$ and ${\cal Z}_\epsilon=\X\setminus {\cal Y}_\epsilon$. Let ${\cal W}_\epsilon=T\setminus D_T$, gluing ${\cal W}_\epsilon$ and ${\cal Z}_\epsilon$ along the boundary $\partial D_T=\Gamma$ one gets the symmetric translation surface $S$ of genus two. (See Figure 1,  the opposite sides of all the parallelograms there are identified.)

By virtue of Proposition \ref{BFK}, one has
\begin{equation}\label{eqN1}
{\rm det}^*\Delta_{{\cal X}}=\frac{{\rm Area}(\X)}{{\rm length}(\Gamma)}{\rm det}(\Delta, {\cal Y}_\epsilon){\rm det}(\Delta, {\cal Z}_\epsilon)
{\rm det}^*({\cal N}_{{\cal Y}_\epsilon}+{\cal N}_{{\cal Z}_\epsilon})\,
\end{equation}
\begin{equation}\label{eqN2}
{\rm det}^*\Delta_S=\frac{2{\rm Area}(T)}{{\rm length} (\Gamma)}{\rm det}(\Delta, {\cal W}_\epsilon){\rm det}(\Delta, {\cal Z}_\epsilon){\rm det}^*
({\cal N}_{{\cal W}_\epsilon}+{\cal N}_{{\cal Z}_\epsilon})
\end{equation}
and, therefore,
\begin{equation}\label{eqN3}
{\rm det}^*\Delta_\X=\frac{{\rm Area}(\X)}{2{\rm Area}(T)}\frac{{\rm det}^*({\cal N}_{{\cal Y}_\epsilon}+{\cal N}_{{\cal Z}_\epsilon})}{{\rm det}^*
({\cal N}_{{\cal W}_\epsilon}+{\cal N}_{{\cal Z}_\epsilon})}\frac{{\rm det}(\Delta, {\cal Y}_\epsilon)}{{\rm det}(\Delta, {\cal W}_\epsilon)}{\rm det}^*\Delta_S\,.
\end{equation}

Using Lemma \ref{Wen}, we have
$$\log \frac{{\rm det}^*({\cal N}_{{\cal Y}_\epsilon}+{\cal N}_{{\cal Z}_\epsilon})}{{\rm det}^*
({\cal N}_{{\cal W}_\epsilon}+{\cal N}_{{\cal Z}_\epsilon})}=\log\frac{{\rm det}^*(|\nu|+{\bf O}_1(\epsilon)+{\cal N}_{{\cal Z}_\epsilon})}{{\rm det}^*
(|\nu|+{\bf O}_2(\epsilon)+{\cal N}_{{\cal Z}_\epsilon})}=$$$$=\int_0^1\frac{d}{dt}\log {\rm det}^*(|\nu|+t{\bf O}_1(\epsilon)+(1-t){\bf O}_2(\epsilon)+{\cal N}_{{\cal Z}_\epsilon})dt=$$$$=\int_0^1{\rm tr}\left[(|\nu|+t{\bf O}_1(\epsilon)+(1-t){\bf O}_2(\epsilon)+{\cal N}_{{\cal Z}_\epsilon})_{\{1\}^\bot}^{-1}({\bf O}_1(\epsilon)-{\bf O}_2(\epsilon))\right]\,dt=O(\epsilon)$$
and, therefore,
$$\frac{{\rm det}^*({\cal N}_{{\cal Y}_\epsilon}+{\cal N}_{{\cal Z}_\epsilon})}{{\rm det}^*
({\cal N}_{{\cal W}_\epsilon}+{\cal N}_{{\cal Z}_\epsilon})}\sim 1$$
as $\epsilon \to 0$.
Due to (\ref{vyrezdisk}) we get
$$ \frac{{\rm det}(\Delta, {\cal Y}_\epsilon)}{{\rm det}(\Delta, {\cal W}_\epsilon)}\sim \frac{{\rm det}^*\Delta_{\cal Y}}{{\rm det}^*\Delta_T}\frac{{\rm Area}(T)}{{\rm Area}(\cal Y)}\,.$$
Finally, from (\ref{AuxAs}) it follows that
$${\rm det}^*\Delta_S\sim\frac{2\kappa_0}{{\rm Area}(T)}\{{\rm det}^*\Delta_T\}^2\epsilon^{\frac{1}{2}}\,.$$
Thus, we conclude from (\ref{eqN3}) that
\begin{equation}\label{fina}
{\rm det}^*\Delta_\X\sim\frac {{\rm Area}(\X)} {{\rm Area}({\cal Y}){\rm Area}(T)}\kappa_0{\rm det}^*\Delta_{\cal Y}{\rm det}^*\Delta_T\epsilon^{\frac{1}{2}}
\end{equation}
as $\epsilon\to 0$.
Comparing (\ref{fina}) and (\ref{Main1}) (with $s:=\epsilon/2$), one arrives at
the main result of the present paper.
\begin{theorem}
The following expression for the constant $\delta_g$ from (\ref{main}) holds true:
\begin{equation}
\delta_{g}=(2\sqrt{2}\kappa_0)^{g-1}\delta_1^{g},
\end{equation}
where $\kappa_0$ is given by (\ref{kappa0}) and $\delta_1$  is from (\ref{delta1}).
\end{theorem}

\end{document}